# Full-waveform inversion in three-dimensional PML-truncated elastic media

Arash Fathi[a,b], Loukas F. Kallivokas[a,b,*], Babak Poursartip[a]

[a]*Department of Civil, Architectural and Environmental Engineering*
[b]*The Institute for Computational Engineering & Sciences*
*The University of Texas at Austin, Austin, TX, USA*

## Abstract

We are concerned with high-fidelity subsurface imaging of the soil, which commonly arises in geotechnical site characterization and geophysical explorations. Specifically, we attempt to image the spatial distribution of the Lamé parameters in semi-infinite, three-dimensional, arbitrarily heterogeneous formations, using surficial measurements of the soil's response to probing elastic waves. We use the complete waveform response of the medium to derive the inverse problem, by using a partial-differential-equation (PDE)-constrained optimization approach, directly in the time-domain, to minimize the misfit between the observed response of the medium at select measurement locations, and a computed response corresponding to a trial distribution of the Lamé parameters. We discuss strategies that lend algorithmic robustness to our proposed inversion scheme. To limit the computational domain to the size of interest, we employ perfectly-matched-layers (PMLs).

In order to resolve the forward problem, we use a recently developed hybrid finite element approach, where a displacement-stress formulation for the PML is coupled to a standard displacement-only formulation for the interior domain, thus leading to a computationally cost-efficient scheme. Time-integration is accomplished by using an explicit Runge-Kutta scheme, which is well-suited for large-scale problems on parallel computers.

We verify the accuracy of the material gradients obtained via our proposed scheme, and report numerical results demonstrating successful reconstruction of the two Lamé parameters for both smooth and sharp profiles.

[*]Corresponding author.
*Email addresses:* `arash.fathi@utexas.edu` (Arash Fathi), `loukas@mail.utexas.edu` (Loukas F. Kallivokas), `babakp@utexas.edu` (Babak Poursartip)

## 1. Introduction

Seismic inversion refers to the process of identification of material properties in geological formations [8, 31, 36]. The problem arises predominantly in exploration geophysics [14, 23, 25, 33] and geotechnical site characterization [20]; it belongs to the broader class of inverse medium problems: waves, whether of acoustic, elastic, or electromagnetic nature, are used to interrogate a medium, and the medium's response to the probing is subsequently used to image the spatial distribution of properties (e.g., Lamé parameters, or wave velocities) [4, 15, 21]. Mathematically, algorithmically, and computationally, inverse medium problems are challenging, especially, when no *a priori* constraining assumption is made on the spatial variability of the medium's properties. The challenges are further compounded when the underlying physics is time-dependent, and involves more than a single distributed parameter to be inverted for, as in seismic inversion.

Due to the complexity of the inverse problem at hand, most techniques to date rely on simplifying assumptions, aiming at rendering a solution to the problem more tractable. These assumptions can be divided into five categories: a) assumptions regarding the dimensionality of the problem, whereby the original problem is reduced to a two-dimensional [17, 20, 22], or a one-dimensional problem [26]; b) assuming that the dominant portion of the wave energy on the ground surface is transported through Rayleigh waves, and thus, disregarding other wave types, such as compressional and shear waves, as is the case in the Spectral-Analysis-of-Surface-Waves (SASW) and its variants (MASW) [35]; c) inverting for only one parameter, as is done in [1, 11, 28, 29], where inversion was attempted only for the shear wave velocity, assuming the compressional wave velocity (or an equivalent counterpart) is known; d) assumptions concerning the truncation boundaries, which are oftentimes, grossly simplified due to the complexity associated with the rigorous treatment of these boundaries [40]; and e) idealizing the soil body, which is a porous and lossy medium, as an elastic solid and neglecting its attenuative properties[1] [12]. Over the past decade, continued advances in both algorithms and computer architectures have allowed the gradual removal of the limitations of existing methodologies. However, a robust methodology, especially for the time-dependent elastic case remains, by and large, elusive.

Among the recent works on inversion, which are similar in character to ours, we refer to Pratt et al. [32] who considered two-dimensional acoustic inversion in the frequency domain, and Epanomeritakis et al. [11] where full-waveform inversion has attempted for three-dimensional time-domain elastodynamics, where a simple boundary condition was used for domain truncation. Kang and Kallivokas [21] considered the problem for the two-dimensional time-domain acoustic case, and used PMLs to accurately account for domain truncation. Kucukcoban [22] extended the work of Kang and Kallivokas to two-dimensional elastodynamics, and reported successful reconstruction of the two Lamé parameters for models involving synthetic data. Bramwell [7] used a discontinuous Petrov-Galerkin (DPG) method in the frequency domain, endowed with PMLs, for seismic tomography problems, advocating the DPG scheme over conventional continuous Galerkin methods, since it results in less numerical pollution. Recently, Jung et al. [18, 19] used an extended finite element method (XFEM) to explicitly parameterize the boundaries of scatterers for two-dimensional problems in elastodynamics. Their approach seems promising especially for the identification of voids.

---

[1]See [2] for a full-waveform-inversion-based approach, using a generalized Maxwell model for lossy soils in a one-dimensional setting.

In this paper, we discuss a systematic framework for the numerical resolution of the inverse medium problem, directly in the time-domain, in the context of geotechnical site characterization. The goal is to image the arbitrarily heterogeneous material profile of a probed soil medium, using complete waveforms[2] of its response to interrogating elastic waves, originating from the ground surface. To this end, the response of the soil medium to active sources (Vibroseis equipment) is collected by receivers (geophones) dispersed over the formation's surface, as shown in Figure 1(a). Arriving at a material profile is then accomplished by minimizing the difference between the collected response at receiver locations, and a computed response corresponding to a trial distribution of the material parameters. Imaging near-surface deposits brings additional difficulties, typically not encountered in exploration geophysics. In geophysical explorations, the probing is over large length scales; thus, an accurate domain termination tool may not play a critical role. However, in geotechnical site characterization, one, typically, deals with a much smaller domain. Moreover, obtaining a high-fidelity image of the near-surface deposits has practical significance for the safe founding of infrastructure components; thus, having accurate termination conditions becomes critical. In this vein, and in the presence of heterogeneity, using PMLs for domain termination is the best available option, and is thus used in this work. The PML is a buffer zone that surrounds the domain of interest, and enforces the decay of outgoing waves. Figure 1(b) shows a computational model, obtained through the introduction of PMLs at the truncation boundaries.

In order to address all the difficulties outlined earlier, we integrate recent advances in several areas. Specifically, we use (a) a recently developed state-of-the-art parallel wave simulation tool for domains terminated by PMLs, which renders the computational model associated with the near surface geotechnical investigations finite [13]; (b) a partial-differential-equation (PDE)-constrained optimization framework through which the minimization of the difference between the collected response at receiver locations and a computed response corresponding to a trial distribution of the material properties is achieved [30]; (c) regularization schemes to alleviate the ill-posedness inherent in inverse problems; (d) continuation schemes that lend algorithmic robustness [21]; and (e) a biasing scheme that accelerates the convergence of the $\lambda$-profile for robust simultaneous inversion of both Lamé parameters

---

[2]Using the complete waveform (complete recorded response) results in a full-waveform inversion approach.

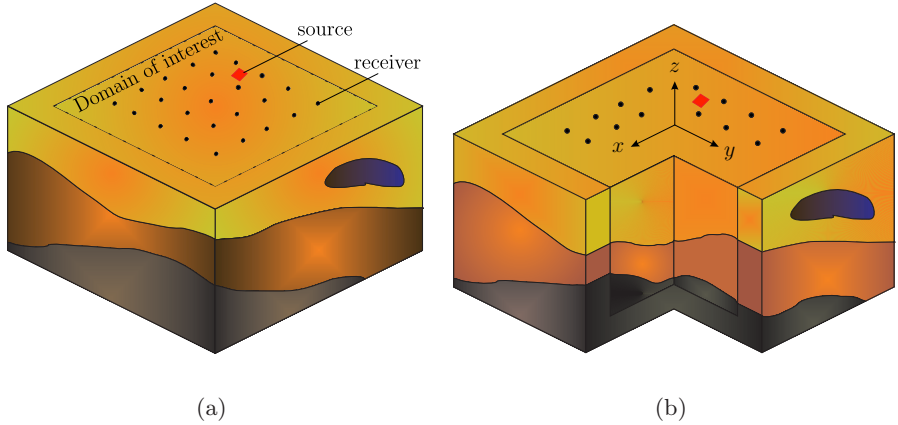

(a)                  (b)

Figure 1: Problem definition: (a) interrogation of a heterogeneous semi-infinite domain by an active source; and (b) computational model truncated from the semi-infinite medium via the introduction of PMLs.

[22].

The remainder of this paper is organized as follows: first, we discuss the numerical resolution of the forward problem. Next, we discuss the mathematical and numerical aspects of the underlying inverse medium problem, where we derive the adjoint and control problems, and discuss strategies that invite robustness. We then verify the accuracy of the material gradients that we compute, by comparing them with directional finite differences. We report on numerical experiments, using synthetic data, targeting the reconstruction of both smooth and sharp profiles. Lastly, we conclude with summary remarks.

## 2. Forward wave simulation in a 3D PML-truncated medium

In the forward problem, we are concerned with the propagation of elastic waves in a three-dimensional, semi-infinite, arbitrarily heterogeneous elastic medium. To keep the computation feasible, one needs to limit the extent of the computational domain. This can be accomplished by placing PMLs at the truncation boundaries. Theoretically, the truncation boundary is reflectionless, and when outgoing waves pass through the interface, they get attenuated within the PML buffer zone. This concept is schematically captured in Figure 2.

For the numerical resolution of the forward problem, we use a recently developed hybrid approach that uses a displacement-stress formulation for the PML buffer, coupled with a standard displacement-only formulation in the interior domain[3]. This approach results in a computationally cost-efficient scheme, due to the treatment of the interior domain with a standard displacement-only elastodynamics formulation. We refer to [13] and references therein for the complete development and parallel implementation of the method. Herein, we only outline the resulting coupled system of equations.

---

[3]The terms "interior domain" and "regular domain" are used interchangeably throughout this article.

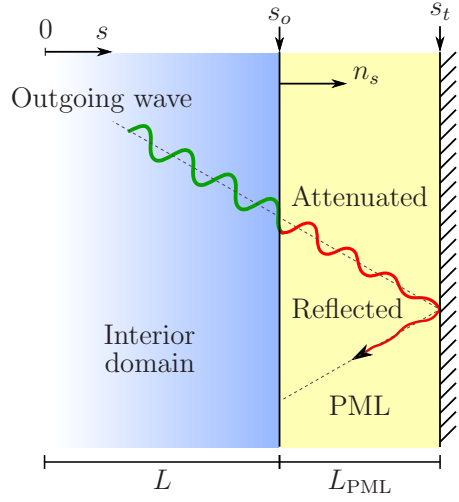

Figure 2: A PML truncation boundary in the direction of coordinate $s$.

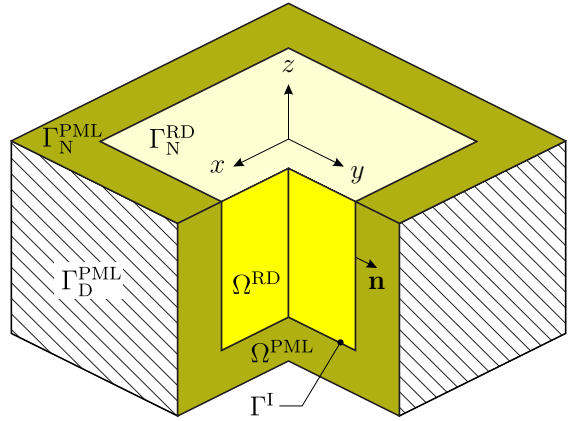

Figure 3: PML-truncated semi-infinite domain.

Accordingly, find $\mathbf{u}(\mathbf{x},t)$ in $\Omega^{\mathrm{RD}} \cup \Omega^{\mathrm{PML}}$, and $\mathbf{S}(\mathbf{x},t)$ in $\Omega^{\mathrm{PML}}$ (see Figure 3 for domain and boundary designations), where $\mathbf{u}$ and $\mathbf{S}$ reside in appropriate function spaces and:

$$\mathbf{div}\left\{\mu\left[\nabla\mathbf{u}+(\nabla\mathbf{u})^T\right]+\lambda(\operatorname{div}\mathbf{u})\mathcal{I}\right\}+\mathbf{b}=\rho\ddot{\mathbf{u}} \qquad \text{in } \Omega^{\text{RD}}\times\mathsf{J}, \tag{1a}$$

$$\mathbf{div}\left(\dot{\mathbf{S}}^T\Lambda_e+\mathbf{S}^T\Lambda_p+\bar{\mathbf{S}}^T\Lambda_w\right)=\rho\left(a\ddot{\mathbf{u}}+b\dot{\mathbf{u}}+c\mathbf{u}+d\bar{\mathbf{u}}\right) \qquad \text{in } \Omega^{\text{PML}}\times\mathsf{J}, \tag{1b}$$

$$a\ddot{\mathbf{S}}+b\dot{\mathbf{S}}+c\mathbf{S}+d\bar{\mathbf{S}}=$$
$$\mu\left[(\nabla\dot{\mathbf{u}})\Lambda_e+\Lambda_e(\nabla\dot{\mathbf{u}})^T+(\nabla\mathbf{u})\Lambda_p+\Lambda_p(\nabla\mathbf{u})^T+(\nabla\bar{\mathbf{u}})\Lambda_w+\Lambda_w(\nabla\bar{\mathbf{u}})^T\right]+$$
$$\lambda\left[\operatorname{div}(\Lambda_e\dot{\mathbf{u}})+\operatorname{div}(\Lambda_p\mathbf{u})+\operatorname{div}(\Lambda_w\bar{\mathbf{u}})\right]\mathcal{I} \qquad \text{in } \Omega^{\text{PML}}\times\mathsf{J}. \tag{1c}$$

The system is initially at rest, and subject to the following boundary and interface conditions:

$$\left\{\mu\left[\nabla\mathbf{u}+(\nabla\mathbf{u})^T\right]+\lambda(\operatorname{div}\mathbf{u})\mathcal{I}\right\}\mathbf{n}=\boldsymbol{g}_n \qquad \text{on } \Gamma_{\text{N}}^{\text{RD}}\times\mathsf{J}, \tag{2a}$$

$$(\dot{\mathbf{S}}^T\Lambda_e+\mathbf{S}^T\Lambda_p+\bar{\mathbf{S}}^T\Lambda_w)\mathbf{n}=\mathbf{0} \qquad \text{on } \Gamma_{\text{N}}^{\text{PML}}\times\mathsf{J}, \tag{2b}$$

$$\mathbf{u}=\mathbf{0} \qquad \text{on } \Gamma_{\text{D}}^{\text{PML}}\times\mathsf{J}, \tag{2c}$$

$$\mathbf{u}^{\text{RD}}=\mathbf{u}^{\text{PML}} \qquad \text{on } \Gamma^{\text{I}}\times\mathsf{J}, \tag{2d}$$

$$\left\{\mu\left[\nabla\mathbf{u}+(\nabla\mathbf{u})^T\right]+\lambda(\operatorname{div}\mathbf{u})\mathcal{I}\right\}\mathbf{n}=(\dot{\mathbf{S}}^T\Lambda_e+\mathbf{S}^T\Lambda_p+\bar{\mathbf{S}}^T\Lambda_w)\mathbf{n} \qquad \text{on } \Gamma^{\text{I}}\times\mathsf{J}, \tag{2e}$$

where temporal and spatial dependencies are suppressed for brevity; $\mathbf{u}$ is the displacement vector, $\rho$ denotes mass density of the medium, $\lambda$ and $\mu$ are the two Lamé parameters, $\mathcal{I}$ is the second-order identity tensor, $\dot{\mathbf{S}}$ represents the Cauchy stress tensor, a dot (˙) denotes differentiation with respect to time, and a bar (¯) indicates history of the subtended variable[4]; $\Omega^{\text{RD}}$ denotes the interior (regular) domain, $\Omega^{\text{PML}}$ represents the region occupied by the PML buffer zone, $\Gamma^{\text{I}}$ is the interface boundary between the interior and PML domains, $\Gamma_{\text{N}}^{\text{RD}}$ and $\Gamma_{\text{N}}^{\text{PML}}$ denote the free (top surface) boundary of the interior domain and PML, respectively, $\mathsf{J}=(0,T]$ is the time interval of interest, $\mathbf{b}$ denotes body force per unit volume, and $\boldsymbol{g}_n$ is the prescribed surface traction. Moreover, $\Lambda_e$, $\Lambda_p$, and $\Lambda_w$ are the so-called stretch tensors, which enforce dissipation of waves in $\Omega^{\text{PML}}$, and $a$, $b$, $c$, and $d$ are products of certain elements of the stretch tensors [13]. Eq. (1a) is the governing PDE for the interior elastodynamic problem, whereas Eqs. (1b) and (1c) are the equilibrium and combined kinematic and constitutive equations, respectively, for the PML buffer.

Next, we seek a weak solution, corresponding to the strong form of (1) and (2), in the Galerkin sense. Specifically, we take the inner products of (1a) and (1b) with (vector) test function $\tilde{\mathbf{w}}(\mathbf{x})$, and integrate by parts over their corresponding domains. Incorporating (2d-2e) eliminates the interface boundary terms and results in (3a). Next, we take the inner product of (1c) with (tensor) test function $\tilde{\mathbf{T}}(\mathbf{x})$; there results (3b). We refer to [13] for details.

---

[4]For instance, $\bar{\mathbf{u}}(\mathbf{x},t)=\int_0^t\mathbf{u}(\mathbf{x},\tau)\,d\tau$.

Accordingly, find $\mathbf{u} \in \mathbf{H}^1(\Omega) \times \mathsf{J}$, and $\mathbf{S} \in \mathcal{L}^2(\Omega) \times \mathsf{J}$, such that:

$$\int_{\Omega^{\text{RD}}} \nabla \tilde{\mathbf{w}} : \left\{ \mu \left[ \nabla \mathbf{u} + (\nabla \mathbf{u})^T \right] + \lambda (\operatorname{div} \mathbf{u}) \mathcal{I} \right\} d\Omega + \int_{\Omega^{\text{PML}}} \nabla \tilde{\mathbf{w}} : \left( \dot{\mathbf{S}}^T \Lambda_e + \mathbf{S}^T \Lambda_p + \bar{\mathbf{S}}^T \Lambda_w \right) d\Omega$$

$$+ \int_{\Omega^{\text{RD}}} \tilde{\mathbf{w}} \cdot \rho \ddot{\mathbf{u}} \, d\Omega + \int_{\Omega^{\text{PML}}} \tilde{\mathbf{w}} \cdot \rho \left( a\ddot{\mathbf{u}} + b\dot{\mathbf{u}} + c\mathbf{u} + d\bar{\mathbf{u}} \right) d\Omega = \int_{\Gamma_N^{\text{RD}}} \tilde{\mathbf{w}} \cdot \boldsymbol{g}_n \, d\Gamma + \int_{\Omega^{\text{RD}}} \tilde{\mathbf{w}} \cdot \mathbf{b} \, d\Omega, \quad (3a)$$

$$\int_{\Omega^{\text{PML}}} \tilde{\mathbf{T}} : \left( a\ddot{\mathbf{S}} + b\dot{\mathbf{S}} + c\mathbf{S} + d\bar{\mathbf{S}} \right) d\Omega$$

$$= \int_{\Omega^{\text{PML}}} \tilde{\mathbf{T}} : \mu \left[ (\nabla \dot{\mathbf{u}}) \Lambda_e + \Lambda_e (\nabla \dot{\mathbf{u}})^T + (\nabla \mathbf{u}) \Lambda_p + \Lambda_p (\nabla \mathbf{u})^T + (\nabla \bar{\mathbf{u}}) \Lambda_w + \Lambda_w (\nabla \bar{\mathbf{u}})^T \right]$$

$$+ \tilde{\mathbf{T}} : \lambda \left[ \operatorname{div}(\Lambda_e \dot{\mathbf{u}}) + \operatorname{div}(\Lambda_p \mathbf{u}) + \operatorname{div}(\Lambda_w \bar{\mathbf{u}}) \right] \mathcal{I} \, d\Omega, \quad (3b)$$

for every $\tilde{\mathbf{w}} \in \mathbf{H}^1(\Omega)$ and $\tilde{\mathbf{T}} \in \mathcal{L}^2(\Omega)$, where $\boldsymbol{g}_n \in \mathbf{L}^2(\Omega) \times \mathsf{J}$, and $\mathbf{b} \in \mathbf{L}^2(\Omega) \times \mathsf{J}$. Function spaces for scalar-valued ($v$), vector-valued ($\mathbf{v}$), and tensor-valued ($\mathcal{A}$) functions are defined as:

$$L^2(\Omega) = \left\{ v : \int_\Omega |v|^2 d\mathbf{x} < \infty \right\}, \quad (4a)$$

$$\mathbf{L}^2(\Omega) = \left\{ \mathbf{v} : \mathbf{v} \in (L^2(\Omega))^3 \right\}, \quad (4b)$$

$$\mathcal{L}^2(\Omega) = \left\{ \mathcal{A} : \mathcal{A} \in (L^2(\Omega))^{3\times 3} \right\}, \quad (4c)$$

$$H^1(\Omega) = \left\{ v : \int_\Omega \left( |v|^2 + |\nabla v|^2 \right) d\mathbf{x} < \infty, \ v|_{\Gamma_D^{\text{PML}}} = 0 \right\}, \quad (4d)$$

$$\mathbf{H}^1(\Omega) = \left\{ \mathbf{v} : \mathbf{v} \in (H^1(\Omega))^3 \right\}. \quad (4e)$$

In order to resolve (3) numerically, we use standard finite-dimensional subspaces. Specifically, we introduce finite-dimensional subspaces $\Xi_h \subset \mathbf{H}^1(\Omega)$ and $\Upsilon_h \subset \mathcal{L}^2(\Omega)$, with basis $\boldsymbol{\Phi} = (\Phi_1, \Phi_2, \ldots, \Phi_N)^T$ and $\boldsymbol{\Psi} = (\Psi_1, \Psi_2, \ldots, \Psi_M)^T$, respectively, where $N$ and $M$ denote the dimension of of the corresponding vector space. We then approximate $\mathbf{u}(\mathbf{x}, t)$ with $\mathbf{u}_h(\mathbf{x}, t) \in \Xi_h \times \mathsf{J}$, and $\mathbf{S}(\mathbf{x}, t)$ with $\mathbf{S}_h(\mathbf{x}, t) \in \Upsilon_h \times \mathsf{J}$. In a similar fashion, we approximate the test functions, $\tilde{\mathbf{w}}(\mathbf{x})$ with $\tilde{\mathbf{w}}_h(\mathbf{x}) \in \Xi_h$, and $\tilde{\mathbf{T}}(\mathbf{x})$ with $\tilde{\mathbf{T}}_h(\mathbf{x}) \in \Upsilon_h$. It can be shown [13] that the following semi-discrete form results:

$$\mathbf{M}\ddot{\mathbf{d}}^{\mathrm{st}} + \mathbf{C}\dot{\mathbf{d}}^{\mathrm{st}} + \mathbf{K}\mathbf{d}^{\mathrm{st}} + \mathbf{G}\bar{\mathbf{d}}^{\mathrm{st}} = \mathbf{f}^{\mathrm{st}}, \tag{5a}$$

$$\bar{\mathbf{d}}^{\mathrm{st}} = \int_0^t \mathbf{d}^{\mathrm{st}}(\tau)|_{\mathrm{PML}} \, d\tau, \tag{5b}$$

where spatial and temporal dependencies are suppressed for brevity; $\mathbf{M}$, $\mathbf{C}$, $\mathbf{K}$, $\mathbf{G}$, are system matrices, $\mathbf{d}^{\mathrm{st}} = (\mathbf{u}_h^T, \mathbf{S}_h^T)^T$ is the vector of nodal unknowns comprising displacements in $\bar{\Omega}^{\mathrm{RD}} \cup \bar{\Omega}^{\mathrm{PML}}$, and stress components only in $\bar{\Omega}^{\mathrm{PML}}$, and $\mathbf{f}^{\mathrm{st}}$ is the vector of applied forces.

The matrix $\mathbf{M}$ has a block-diagonal structure, and can be diagonalized if one employs spectral elements with a Legendre-Gauss-Lobatto (LGL) quadrature rule, which then enables explicit time integration of (5). In this regard, we express (5) as a first-order system:

$$\frac{d}{dt}\begin{bmatrix} \mathbf{x}_0 \\ \mathbf{x}_1 \\ \mathbf{M}\mathbf{x}_2 \end{bmatrix} = \begin{bmatrix} \mathbf{0} & \mathbf{I} & \mathbf{0} \\ \mathbf{0} & \mathbf{0} & \mathbf{I} \\ -\mathbf{G} & -\mathbf{K} & -\mathbf{C} \end{bmatrix} \begin{bmatrix} \mathbf{x}_0 \\ \mathbf{x}_1 \\ \mathbf{x}_2 \end{bmatrix} + \begin{bmatrix} \mathbf{0} \\ \mathbf{0} \\ \mathbf{f}^{\mathrm{st}} \end{bmatrix}, \tag{6}$$

where $\mathbf{x}_0 = \bar{\mathbf{d}}^{\mathrm{st}}$, $\mathbf{x}_1 = \mathbf{d}^{\mathrm{st}}$, and $\mathbf{x}_2 = \dot{\mathbf{d}}^{\mathrm{st}}$. We then use an explicit fourth-order Runge-Kutta (RK-4) method for integrating (6) in time.

## 3. The inverse medium problem

Our goal is to find the distribution of the Lamé parameters $\lambda(\mathbf{x})$ and $\mu(\mathbf{x})$ within the elastic soil medium. We consider sources, and the response recorded at receivers on the ground surface, as known. The inverse medium problem can thus be formulated as the minimization of the difference (or misfit) between the measured response at receiver locations, and a computed response corresponding to a trial distribution of the material parameters. The misfit minimization should honor the physics of the problem, which is idealized by the forward problem, stated in the preceding section. Mathematically, this can be cast as a PDE-constrained optimization problem:

$$\min_{\lambda,\mu} \; \mathcal{J}(\lambda,\mu) := \frac{1}{2} \sum_{j=1}^{N_r} \int_0^T \int_{\Gamma_m} (\mathbf{u} - \mathbf{u}_m) \cdot (\mathbf{u} - \mathbf{u}_m) \, \delta(\mathbf{x} - \mathbf{x}_j) \, d\Gamma \, dt + \mathcal{R}(\lambda,\mu), \tag{7}$$

where $\mathbf{u}$ is the solution of the forward problem governed by the initial- and boundary-value problem (1), (2).

In the above, $\mathcal{J}$ is the objective functional[5], $N_r$ denotes the total number of receivers, $T$ is the total observation time, $\Gamma_m$ is the part of the ground surface where the receiver response, $\mathbf{u}_m$, has been recorded, $\delta(\mathbf{x} - \mathbf{x}_j)$ is the Dirac delta function, which enables measurements at receiver locations $\mathbf{x}_j$, and $\mathcal{R}(\lambda, \mu)$ is the regularization term, which is discussed below.

---
[5]We use $\mathbf{J}$ to indicate the corresponding discrete objective functional. See [6, 10] for other possibilities.

Inverse problems suffer from solution multiplicity, which, in general, is due to the presence of insufficient data. This makes the problem ill-posed in the Hadamard sense. Regularization of the solution by using the Tikhonov (TN) [37], or, the Total Variation (TV) [34] scheme are among common strategies to alleviate ill-posedness. The Tikhonov regularization, denoted by $\mathcal{R}^{TN}(\lambda, \mu)$, penalizes large material gradients and, thus, precludes spatially rapid material variations from becoming solutions to the inverse medium problem. It is defined as:

$$\mathcal{R}^{TN}(\lambda, \mu) = \frac{R_\lambda}{2} \int_{\Omega^{\text{RD}}} \nabla\lambda \cdot \nabla\lambda \, \mathrm{d}\Omega + \frac{R_\mu}{2} \int_{\Omega^{\text{RD}}} \nabla\mu \cdot \nabla\mu \, \mathrm{d}\Omega, \tag{8}$$

where $R_\lambda$ and $R_\mu$ are the so-called $\lambda$- and $\mu$-regularization factor, respectively, and control the amount of penalty imposed via (8) on the gradients of $\lambda$ and $\mu$. By construction, TN regularization results in material reconstructions with smooth variations. Consequently, sharp interfaces may not be captured well when using the TN scheme. The TV regularization, however, works better for imaging profiles involving sharp interfaces, as it typically preserves edges. It is defined as:

$$\mathcal{R}^{TV}(\lambda, \mu) = \frac{R_\lambda}{2} \int_{\Omega^{\text{RD}}} (\nabla\lambda \cdot \nabla\lambda + \epsilon)^{\frac{1}{2}} \, \mathrm{d}\Omega + \frac{R_\mu}{2} \int_{\Omega^{\text{RD}}} (\nabla\mu \cdot \nabla\mu + \epsilon)^{\frac{1}{2}} \, \mathrm{d}\Omega, \tag{9}$$

where the parameter $\epsilon$ makes $\mathcal{R}^{TV}$ differentiable when either $\nabla\lambda \cdot \nabla\lambda = 0$, or, $\nabla\mu \cdot \nabla\mu = 0$.

For computing the first-order optimality conditions for (7), we use the (formal) Lagrangian approach [41] to impose the PDE-constraint in its weak form. These are necessary conditions that must be satisfied at a solution of (7). Specifically, we introduce Lagrange multiplier vector function $\mathbf{w} \in \mathbf{H}^1(\Omega)$, and Lagrange multiplier tensor function $\mathbf{T} \in \mathcal{L}^2(\Omega)$ to enforce the initial- and boundary-value problem (1), (2), which admits the weak form (3). The Lagrangian functional becomes:

$$\mathcal{L}(\mathbf{u},\mathbf{S},\mathbf{w},\mathbf{T},\lambda,\mu) = \frac{1}{2}\sum_{j=1}^{N_r}\int_0^T\int_{\Gamma_m}(\mathbf{u}-\mathbf{u}_m)\cdot(\mathbf{u}-\mathbf{u}_m)\,\delta(\mathbf{x}-\mathbf{x}_j)\,\mathrm{d}\Gamma\,dt + \mathcal{R}(\lambda,\mu)$$

$$-\int_0^T\int_{\Omega^{\mathrm{RD}}}\nabla\mathbf{w}:\{\mu\left[\nabla\mathbf{u}+(\nabla\mathbf{u})^T\right]+\lambda(\operatorname{div}\mathbf{u})\mathcal{I}\}\,\mathrm{d}\Omega\,dt$$

$$-\int_0^T\int_{\Omega^{\mathrm{PML}}}\nabla\mathbf{w}:\left(\dot{\mathbf{S}}^T\Lambda_e+\mathbf{S}^T\Lambda_p+\bar{\mathbf{S}}^T\Lambda_w\right)\,\mathrm{d}\Omega\,dt - \int_0^T\int_{\Omega^{\mathrm{RD}}}\mathbf{w}\cdot\rho\ddot{\mathbf{u}}\,\mathrm{d}\Omega\,dt$$

$$-\int_0^T\int_{\Omega^{\mathrm{PML}}}\mathbf{w}\cdot\rho\left(a\ddot{\mathbf{u}}+b\dot{\mathbf{u}}+c\mathbf{u}+d\bar{\mathbf{u}}\right)\,\mathrm{d}\Omega\,dt + \int_0^T\int_{\Gamma_N^{\mathrm{RD}}}\mathbf{w}\cdot\boldsymbol{g}_n\,\mathrm{d}\Gamma\,dt$$

$$+\int_0^T\int_{\Omega^{\mathrm{RD}}}\mathbf{w}\cdot\mathbf{b}\,\mathrm{d}\Omega\,dt - \int_0^T\int_{\Omega^{\mathrm{PML}}}\mathbf{T}:\left(a\ddot{\mathbf{S}}+b\dot{\mathbf{S}}+c\mathbf{S}+d\bar{\mathbf{S}}\right)\,\mathrm{d}\Omega\,dt$$

$$+\int_0^T\int_{\Omega^{\mathrm{PML}}}\mathbf{T}:\mu\left[(\nabla\dot{\mathbf{u}})\Lambda_e+\Lambda_e(\nabla\dot{\mathbf{u}})^T+(\nabla\mathbf{u})\Lambda_p+\Lambda_p(\nabla\mathbf{u})^T+(\nabla\bar{\mathbf{u}})\Lambda_w+\Lambda_w(\nabla\bar{\mathbf{u}})^T\right]$$

$$+\mathbf{T}:\lambda\left[\operatorname{div}(\Lambda_e\dot{\mathbf{u}})+\operatorname{div}(\Lambda_p\mathbf{u})+\operatorname{div}(\Lambda_w\bar{\mathbf{u}})\right]\mathcal{I}\,\mathrm{d}\Omega\,dt, \tag{10a}$$

where now $\mathbf{u}$, $\mathbf{S}$, $\lambda$, and $\mu$ are treated as independent variables.

*3.1. Optimality system*

We now use the Lagrangian functional (10) as a tool to compute the optimality system for (7). To this end, the Gâteaux derivative[6] (or first variation) of the Lagrangian functional with respect to all variables must vanish. This process is discussed next.

*3.1.1. The state problem*

Taking the derivatives of the Lagrangian functional $\mathcal{L}$ with respect to $\mathbf{w}$ and $\mathbf{T}$ in directions $\tilde{\mathbf{w}}\in\mathbf{H}^1(\Omega)$ and $\tilde{\mathbf{T}}\in\mathcal{L}^2(\Omega)$, and setting it to zero, results in the state problem, which is identical to (3). That is:

$$\mathcal{L}'(\mathbf{u},\mathbf{S},\mathbf{w},\mathbf{T},\lambda,\mu)(\tilde{\mathbf{w}},\tilde{\mathbf{T}}) = 0. \tag{11}$$

*3.1.2. The adjoint problem*

We now take the derivative of $\mathcal{L}$ with respect to $\mathbf{u}$ and $\mathbf{S}$ in directions $\tilde{\mathbf{u}}\in\mathbf{H}^1(\Omega)$ and $\tilde{\mathbf{S}}\in\mathcal{L}^2(\Omega)$. This yields:

---

[6]See Appendix A for the definition and notation.

$$\begin{aligned}
\mathcal{L}'(\mathbf{u},\mathbf{S},\mathbf{w},\mathbf{T},\lambda,\mu)(\tilde{\mathbf{u}},\tilde{\mathbf{S}}) = &\sum_{j=1}^{N_r} \int_0^T \int_{\Gamma_m} \tilde{\mathbf{u}}\cdot(\mathbf{u}-\mathbf{u}_m)\,\delta(\mathbf{x}-\mathbf{x}_j)\,\mathrm{d}\Gamma\,dt \\
&- \int_0^T \int_{\Omega^{\mathrm{RD}}} \nabla\mathbf{w} : \left\{\mu\left[\nabla\tilde{\mathbf{u}} + (\nabla\tilde{\mathbf{u}})^T\right] + \lambda(\mathrm{div}\,\tilde{\mathbf{u}})\mathcal{I}\right\}\,\mathrm{d}\Omega\,dt \\
&- \int_0^T \int_{\Omega^{\mathrm{PML}}} \nabla\mathbf{w} : \left(\dot{\tilde{\mathbf{S}}}^T \Lambda_e + \tilde{\mathbf{S}}^T \Lambda_p + \bar{\tilde{\mathbf{S}}}^T \Lambda_w\right)\,\mathrm{d}\Omega\,dt - \int_0^T \int_{\Omega^{\mathrm{RD}}} \mathbf{w}\cdot\rho\ddot{\tilde{\mathbf{u}}}\,\mathrm{d}\Omega\,dt \\
&- \int_0^T \int_{\Omega^{\mathrm{PML}}} \mathbf{w}\cdot\rho\left(a\ddot{\tilde{\mathbf{u}}} + b\dot{\tilde{\mathbf{u}}} + c\tilde{\mathbf{u}} + d\bar{\tilde{\mathbf{u}}}\right)\,\mathrm{d}\Omega\,dt \\
&- \int_0^T \int_{\Omega^{\mathrm{PML}}} \mathbf{T} : \left(a\ddot{\tilde{\mathbf{S}}} + b\dot{\tilde{\mathbf{S}}} + c\tilde{\mathbf{S}} + d\bar{\tilde{\mathbf{S}}}\right)\,\mathrm{d}\Omega\,dt \\
&+ \int_0^T \int_{\Omega^{\mathrm{PML}}} \mathbf{T} : \mu\left[(\nabla\dot{\tilde{\mathbf{u}}})\Lambda_e + \Lambda_e(\nabla\dot{\tilde{\mathbf{u}}})^T + (\nabla\tilde{\mathbf{u}})\Lambda_p + \Lambda_p(\nabla\tilde{\mathbf{u}})^T + (\nabla\bar{\tilde{\mathbf{u}}})\Lambda_w + \Lambda_w(\nabla\bar{\tilde{\mathbf{u}}})^T\right] \\
&+ \mathbf{T} : \lambda\left[\mathrm{div}(\Lambda_e \dot{\tilde{\mathbf{u}}}) + \mathrm{div}(\Lambda_p \tilde{\mathbf{u}}) + \mathrm{div}(\Lambda_w \bar{\tilde{\mathbf{u}}})\right]\mathcal{I}\,\mathrm{d}\Omega\,dt. \quad (12\mathrm{a})
\end{aligned}$$

Setting the above derivative to zero, and performing integration by parts in time, results in the statement of the weak form of the adjoint problem. That is, find $\mathbf{w} \in \mathbf{H}^1(\Omega) \times \mathsf{J}$, and $\mathbf{T} \in \mathcal{L}^2(\Omega) \times \mathsf{J}$, such that:

$$\int_{\Omega^{\text{RD}}} \nabla \tilde{\mathbf{u}} : \left\{ \mu \left[ \nabla \mathbf{w} + (\nabla \mathbf{w})^T \right] + \lambda (\text{div } \mathbf{w}) \mathcal{I} \right\} \, d\Omega +$$

$$+ \int_{\Omega^{\text{RD}}} \tilde{\mathbf{u}} \cdot \rho \ddot{\mathbf{w}} \, d\Omega + \int_{\Omega^{\text{PML}}} \tilde{\mathbf{u}} \cdot \rho \left( a \ddot{\mathbf{w}} - b \dot{\mathbf{w}} + c \mathbf{w} - d \bar{\mathbf{w}} \right) \, d\Omega$$

$$- \int_{\Omega^{\text{PML}}} \nabla \tilde{\mathbf{u}} : \mu \left[ -\dot{\mathbf{T}} \Lambda_e - \dot{\mathbf{T}}^T \Lambda_e + \mathbf{T} \Lambda_p + \mathbf{T}^T \Lambda_p - \bar{\mathbf{T}} \Lambda_w - \bar{\mathbf{T}}^T \Lambda_w \right]$$

$$+ \lambda \left[ -\dot{\mathbf{T}} : \text{div}(\Lambda_e \tilde{\mathbf{u}}) + \mathbf{T} : \text{div}(\Lambda_p \tilde{\mathbf{u}}) - \bar{\mathbf{T}} : \text{div}(\Lambda_w \tilde{\mathbf{u}}) \right] \mathcal{I} \, d\Omega$$

$$= \sum_{j=1}^{N_r} \int_{\Gamma_m} \tilde{\mathbf{u}} \cdot (\mathbf{u} - \mathbf{u}_m) \, \delta(\mathbf{x} - \mathbf{x}_j) \, d\Gamma, \tag{13a}$$

$$\int_{\Omega^{\text{PML}}} \nabla \dot{\mathbf{w}} : \tilde{\mathbf{S}}^T \Lambda_e - \nabla \mathbf{w} : \tilde{\mathbf{S}}^T \Lambda_p + \nabla \bar{\mathbf{w}} : \tilde{\mathbf{S}}^T \Lambda_w \, d\Omega = \int_{\Omega^{\text{PML}}} \tilde{\mathbf{S}} : \left( a \ddot{\mathbf{T}} - b \dot{\mathbf{T}} + c \mathbf{S} - d \bar{\mathbf{T}} \right) \, d\Omega \tag{13b}$$

for every $\tilde{\mathbf{u}} \in \mathbf{H}^1(\Omega)$ and $\tilde{\mathbf{S}} \in \mathcal{L}^2(\Omega)$, where $\mathbf{w}(\mathbf{x}, T) = \mathbf{0}$, and $\mathbf{T}(\mathbf{x}, T) = \mathbf{0}$.

We remark that the adjoint problem (13) is a final-value problem and, thus, is solved backwards in time[7]; it is driven by the misfit between a computed response, and the measured response at receiver locations. Moreover, the operators implicated in the adjoint problem are very similar to those of the state problem: they involve transposition of the system matrices, and sign reversal for terms involving history, and first-order time derivatives. In this regard, we obtain the following semi-discrete form for the adjoint problem:

$$\mathbf{M} \ddot{\mathbf{d}}^{\text{adj}} - \mathbf{C}^T \dot{\mathbf{d}}^{\text{adj}} + \mathbf{K}^T \mathbf{d}^{\text{adj}} - \mathbf{G}^T \bar{\mathbf{d}}^{\text{adj}} = \mathbf{f}^{\text{adj}}, \tag{14a}$$

$$\bar{\mathbf{d}}^{\text{adj}} = \int_0^t \mathbf{d}^{\text{adj}}(\tau)|_{\text{PML}} \, d\tau, \tag{14b}$$

where superscript "adj" refers to the adjoint problem, $\mathbf{d}^{\text{adj}} = (\mathbf{w}_h^T, \mathbf{T}_h^T)^T$ is the vector of nodal unknowns comprising discrete values of $\mathbf{w}$ in $\bar{\Omega}^{\text{RD}} \cup \bar{\Omega}^{\text{PML}}$ and discrete values of $\mathbf{T}$ only in $\bar{\Omega}^{\text{PML}}$, and $\mathbf{f}^{\text{adj}}$ is a vector comprising the misfit at receiver's locations. Moreover, system matrices $\mathbf{M}$, $\mathbf{C}$, $\mathbf{K}$, $\mathbf{G}$, are identical to those of the forward problem and, thus, with minor adjustments, an implementation of the forward problem can also be used for the solution of the adjoint problem.

---
[7] See [44] for other possibilities, and [24, 38, 39] for alternative approaches.

The matrix $\mathbf{M}$ in (14) can be diagonalized by using spectral elements with a Legendre-Gauss-Lobatto (LGL) quadrature rule, similar to what we did in (6). We rewrite (14) as a first-order system:

$$\frac{d}{dt}\begin{bmatrix}\mathbf{y}_0\\ \mathbf{y}_1\\ \mathbf{M}\mathbf{y}_2\end{bmatrix} = \begin{bmatrix}\mathbf{0} & \mathbf{I} & \mathbf{0}\\ \mathbf{0} & \mathbf{0} & \mathbf{I}\\ \mathbf{G}^T & -\mathbf{K}^T & \mathbf{C}^T\end{bmatrix}\begin{bmatrix}\mathbf{y}_0\\ \mathbf{y}_1\\ \mathbf{y}_2\end{bmatrix} + \begin{bmatrix}\mathbf{0}\\ \mathbf{0}\\ \mathbf{f}^{\mathrm{adj}}\end{bmatrix}, \quad (15)$$

where $\mathbf{y}_0 = \bar{\mathbf{d}}^{\mathrm{adj}}$, $\mathbf{y}_1 = \mathbf{d}^{\mathrm{adj}}$, $\mathbf{y}_2 = \dot{\mathbf{d}}^{\mathrm{adj}}$, with final values $\mathbf{y}_0(T) = \mathbf{0}$, $\mathbf{y}_1(T) = \mathbf{0}$, and $\mathbf{y}_2(T) = \mathbf{0}$. We use an explicit RK-4 method to integrate (15) in time. The scheme is outlined in Appendix B.

*3.1.3. The control problems*

Lastly, we take the derivative of $\mathcal{L}$ with respect to $\lambda$ and $\mu$ in directions $\tilde{\lambda}$ and $\tilde{\mu}$, which yields the reduced gradients with respect to $\lambda$ and $\mu$, respectively. We restrict the reduced gradients to $\Omega^{\mathrm{RD}}$ (The material properties at the interfaces $\Gamma^{\mathrm{I}}$ are extended into the PML buffer.). For the TN regularization, this yields:

$$\mathcal{L}'(\mathbf{u},\mathbf{S},\mathbf{w},\mathbf{T},\lambda,\mu)(\tilde{\lambda}) = R_\lambda \int_{\Omega^{\mathrm{RD}}} \nabla\tilde{\lambda}\cdot\nabla\lambda\,d\Omega - \int_0^T\int_{\Omega^{\mathrm{RD}}} \tilde{\lambda}\,\nabla\mathbf{w}:(\mathrm{div}\,\mathbf{u})\mathcal{I}\,d\Omega\,dt, \quad (16a)$$

$$\mathcal{L}'(\mathbf{u},\mathbf{S},\mathbf{w},\mathbf{T},\lambda,\mu)(\tilde{\mu}) = R_\mu \int_{\Omega^{\mathrm{RD}}} \nabla\tilde{\mu}\cdot\nabla\mu\,d\Omega - \int_0^T\int_{\Omega^{\mathrm{RD}}} \tilde{\mu}\,\nabla\mathbf{u}:\left[\nabla\mathbf{w}+(\nabla\mathbf{w})^T\right]\,d\Omega\,dt. \quad (16b)$$

Setting the above derivatives to zero, results in the control problems. Similarly, for the TV regularization, the control problems read:

$$\mathcal{L}'(\mathbf{u},\mathbf{S},\mathbf{w},\mathbf{T},\lambda,\mu)(\tilde{\lambda}) = R_\lambda \int_{\Omega^{\mathrm{RD}}} \frac{\nabla\tilde{\lambda}\cdot\nabla\lambda}{(\nabla\lambda\cdot\nabla\lambda+\epsilon)^{\frac{1}{2}}}\,d\Omega - \int_0^T\int_{\Omega^{\mathrm{RD}}} \tilde{\lambda}\,\nabla\mathbf{w}:(\mathrm{div}\,\mathbf{u})\mathcal{I}\,d\Omega\,dt, \quad (17a)$$

$$\mathcal{L}'(\mathbf{u},\mathbf{S},\mathbf{w},\mathbf{T},\lambda,\mu)(\tilde{\mu}) = R_\mu \int_{\Omega^{\mathrm{RD}}} \frac{\nabla\tilde{\mu}\cdot\nabla\mu}{(\nabla\mu\cdot\nabla\mu+\epsilon)^{\frac{1}{2}}}\,d\Omega - \int_0^T\int_{\Omega^{\mathrm{RD}}} \tilde{\mu}\,\nabla\mathbf{u}:\left[\nabla\mathbf{w}+(\nabla\mathbf{w})^T\right]\,d\Omega\,dt. \quad (17b)$$

Discretization of either (16) or (17) result in the following form:

$$\tilde{\mathbf{M}}\mathbf{g}^\lambda = R_\lambda \ \mathbf{g}^\lambda_{\text{reg}} + \mathbf{g}^\lambda_{\text{mis}}, \tag{18a}$$

$$\tilde{\mathbf{M}}\mathbf{g}^\mu = R_\mu \ \mathbf{g}^\mu_{\text{reg}} + \mathbf{g}^\mu_{\text{mis}}, \tag{18b}$$

where $\tilde{\mathbf{M}}$ is a mass-like matrix, $\mathbf{g}^\lambda$ and $\mathbf{g}^\mu$ is the vector of discrete values of the (reduced) gradient for $\lambda$ and $\mu$, respectively, and $\mathbf{g}^\lambda_{\text{reg}}$, $\mathbf{g}^\mu_{\text{reg}}$ and $\mathbf{g}^\lambda_{\text{mis}}$, $\mathbf{g}^\mu_{\text{mis}}$ are the associated vectors corresponding to the regularization-part and misfit-part of $\mathbf{g}^\lambda$ and $\mathbf{g}^\mu$. We refer to Appendix C for matrix and vector definitions, and discretization details.

*3.2. The inversion process*

A solution of (7) requires simultaneous satisfaction of the state problem (6), the adjoint problem (15), and the control problems (18). This approach –a full-space method– is, in principle, possible [9]; however, the associated computational cost can be substantial. Alternatively, a reduced-space method may be adopted, in which, discrete material properties are updated iteratively, using a gradient-based minimization scheme. The latter approach is employed here, and is discussed next.

We start with an assumed initial spatial distribution of the control parameters ($\lambda$ and $\mu$), and solve the state problem (6) to obtain $\mathbf{d}^{\text{st}} = (\mathbf{u}_h^T, \mathbf{S}_h^T)^T$. With the misfit known, we solve the adjoint problem (15) and obtain $\mathbf{d}^{\text{adj}} = (\mathbf{w}_h^T, \mathbf{T}_h^T)^T$. With $\mathbf{u}_h$ and $\mathbf{w}_h$ known, the (reduced) material gradients, i.e., $\mathbf{g}^\lambda$ and $\mathbf{g}^\mu$, can be computed from (18). Thus, the vector of material values, at iteration $k+1$, can be computed by using a search direction via:

$$\boldsymbol{\lambda}_{k+1} = \boldsymbol{\lambda}_k + \alpha_k^\lambda \ \boldsymbol{s}_k^\lambda, \tag{19a}$$

$$\boldsymbol{\mu}_{k+1} = \boldsymbol{\mu}_k + \alpha_k^\mu \ \boldsymbol{s}_k^\mu, \tag{19b}$$

where $\boldsymbol{\lambda}$ and $\boldsymbol{\mu}$ comprise the vector of discrete values for $\lambda$ and $\mu$, respectively, $\alpha_k^\lambda$, $\alpha_k^\mu$ are step lengths, and $\boldsymbol{s}_k^\lambda$, $\boldsymbol{s}_k^\mu$ are the search directions for $\boldsymbol{\lambda}_k$ and $\boldsymbol{\mu}_k$. Herein, we use the L-BFGS method to compute the search directions [27][8]. Moreover, to ensure sufficient decrease of the objective functional at each inversion iteration, we employ an Armijo backtracking line search [27], which is outlined in Algorithm 1. The inversion process discussed thus far is summarized in Algorithm 2.

We remark that for the reduced-space method, either (16) or (17) can also be expressed as:

$$\mathcal{L}'(\mathbf{u}, \mathbf{S}, \mathbf{w}, \mathbf{T}, \lambda, \mu)(\tilde{\lambda}) = \mathcal{J}'(\lambda, \mu)(\tilde{\lambda}), \tag{20a}$$

$$\mathcal{L}'(\mathbf{u}, \mathbf{S}, \mathbf{w}, \mathbf{T}, \lambda, \mu)(\tilde{\mu}) = \mathcal{J}'(\lambda, \mu)(\tilde{\mu}), \tag{20b}$$

---

[8]In the numerical experiments that we perform, we store $m = 15$ L-BFGS vectors.

**Algorithm 1** Backtracking line search.
___
1: Choose $\alpha^\lambda$, $\alpha^\mu$, $c_1$, $\rho$      ▷ e.g., $\alpha^\lambda = 1$, $\alpha^\mu = 1$, $c_1 = 10^{-4}$, $\rho = 0.5$
2: **while** $\mathbf{J}(\boldsymbol{\lambda}_k + \alpha^\lambda \, \boldsymbol{s}_k^\lambda, \boldsymbol{\mu}_k + \alpha^\mu \, \boldsymbol{s}_k^\mu) \geq \mathbf{J}(\boldsymbol{\lambda}_k, \boldsymbol{\mu}_k) + c_1(\alpha^\lambda \, \mathbf{g}_k^\lambda \cdot \boldsymbol{s}_k^\lambda + \alpha^\mu \, \mathbf{g}_k^\mu \cdot \boldsymbol{s}_k^\mu)$ **do**
3:      $\alpha^\lambda \leftarrow \rho \alpha^\lambda$
4:      $\alpha^\mu \leftarrow \rho \alpha^\mu$
5: **end while**
6: Terminate with $\alpha_k^\lambda = \alpha^\lambda$, $\alpha_k^\mu = \alpha^\mu$
___

**Algorithm 2** Inversion for Lamé parameters.
___
1: $k \leftarrow 0$
2: Set initial guess for material property vectors $\boldsymbol{\lambda}_k$, $\boldsymbol{\mu}_k$
3: Compute $\mathbf{J}(\boldsymbol{\lambda}_k, \boldsymbol{\mu}_k)$      ▷ Eq. (7)
4: Set convergence tolerance *tol*
5: **while** $\{\mathbf{J}(\boldsymbol{\lambda}_k, \boldsymbol{\mu}_k) > tol\}$ **do**
6:      Solve the state problem for $\mathbf{d}^{\text{st}} = (\mathbf{u}_h^T, \mathbf{S}_h^T)^T$      ▷ Eq. (6)
7:      Solve the adjoint problem for $\mathbf{d}^{\text{adj}} = (\mathbf{w}_h^T, \mathbf{T}_h^T)^T$      ▷ Eq. (15)
8:      Evaluate the discrete reduced gradients $\mathbf{g}_k^\lambda$, $\mathbf{g}_k^\mu$      ▷ Eqs. (18)
9:      Compute search directions $\boldsymbol{s}_k^\lambda$, $\boldsymbol{s}_k^\mu$      ▷ L-BFGS
10:     Choose step lengths $\alpha_k^\lambda$, $\alpha_k^\mu$      ▷ Algorithm 1
11:     Update material property vectors $\boldsymbol{\lambda}_k$, $\boldsymbol{\mu}_k$      ▷ Eq. (19)
12:     $k \leftarrow k + 1$
13: **end while**
___

where the equality in (20) is due to the satisfaction of the state problem. Therefore, the reduced gradients in (18), are, indeed, the gradients of the objective functional with respect to $\lambda$ and $\mu$.

*3.3. Buttressing schemes*

Inverse medium problems are notoriously ill-posed. They suffer from solution multiplicity; that is, material profiles that are very different from each other, and, potentially non-physical, can become solutions to the misfit minimization problem. Regularization of the control parameters alleviates the ill-posedness; however, this alone, may not be adequate when dealing with large-scale complex problems. In this part, we discuss additional strategies that further assist the inversion process, outlined in Algorithm 2, to image complex profiles.

*3.3.1. Regularization factor selection and continuation*

Computation of the (reduced) gradients (18) necessitates selection of the regularization factors $R_\lambda$ and $R_\mu$. A common strategy is to use Morozov's discrepancy principle [42], where a constant value for the regularization factor is used throughout the inversion process. Here, we discuss a simple and practical approach that was initially developed for acoustic inversion [21], and, later, was successfully applied to problems involving elastic inversion [22].

We start by rewriting the discrete control problem (18), either for $\lambda$ or $\mu$, in the following generic form:

$$\tilde{\mathbf{M}} \mathbf{g} = R \, \mathbf{g}_{\text{reg}} + \mathbf{g}_{\text{mis}}, \tag{21}$$

where $\mathbf{g}$ refers to the vector of discrete values of the (reduced) gradient, either for $\lambda$ or $\mu$, $\mathbf{g}_{\text{reg}}$ and $\mathbf{g}_{\text{mis}}$ are the associated vectors corresponding to the regularization-part and misfit-part of $\mathbf{g}$, and $R$ is the regularization factor yet to be determined. The main idea is that the "size" of $R \, \mathbf{g}_{\text{reg}}$ should be proportional to that of $\mathbf{g}_{\text{mis}}$ at each inversion iteration. We define the following unit vectors for the two components of the gradient vector:

$$\mathbf{n}_{\text{reg}} = \frac{\mathbf{g}_{\text{reg}}}{\|\mathbf{g}_{\text{reg}}\|}, \qquad \mathbf{n}_{\text{mis}} = \frac{\mathbf{g}_{\text{mis}}}{\|\mathbf{g}_{\text{mis}}\|}, \tag{22}$$

where $\|\cdot\|$ denotes the Euclidean norm. Equation (21) can then be written as:

$$\tilde{\mathbf{M}} \mathbf{g} = R \, \|\mathbf{g}_{\text{reg}}\| \, \mathbf{n}_{\text{reg}} + \|\mathbf{g}_{\text{mis}}\| \, \mathbf{n}_{\text{mis}} \tag{23a}$$

$$= \|\mathbf{g}_{\text{mis}}\| \left( R \, \frac{\|\mathbf{g}_{\text{reg}}\|}{\|\mathbf{g}_{\text{mis}}\|} \, \mathbf{n}_{\text{reg}} + \mathbf{n}_{\text{mis}} \right)$$

$$= \|\mathbf{g}_{\text{mis}}\| \left( \wp \, \mathbf{n}_{\text{reg}} + \mathbf{n}_{\text{mis}} \right), \tag{23b}$$

where,

$$\wp = R \, \frac{\|\mathbf{g}_{\text{reg}}\|}{\|\mathbf{g}_{\text{mis}}\|}. \tag{23c}$$

In (23b), for the "size" of $\wp \, \mathbf{n}_{\text{reg}}$ to be proportional to $\mathbf{n}_{\text{mis}}$ throughout the entire inversion process, one may choose $0 < \wp \leq 1$. Once a value for $\wp$ has been decided on, $R$ can be computed via:

$$R = \wp \, \frac{\|\mathbf{g}_{\text{mis}}\|}{\|\mathbf{g}_{\text{reg}}\|}, \tag{24}$$

where $\wp$ can take large values (e.g., 0.5)[9] at early stages of inversion and, thus, narrow down the initial search space. As inversion evolves, $\wp$ can be continuously reduced (e.g., down to 0.3) to allow for profile refinement. The suggested values for $\wp$ are based on our experience with various numerical experiments that provide quality solutions for different test cases, and seem to be independent of the dimensionality, size, and discretization of the considered cases.

### 3.3.2. Source-frequency continuation

Using loading sources with low-frequency content result in an overall image of the medium that lacks fine features. To allow for more details, and fine-tune the profile, one needs to use sources with higher frequency content. Thus, the inversion process can be initiated with a signal having a low-frequency content and, then, the frequency range can be increased progressively as inversion evolves. This can be achieved by using a set of probing signals, ordered such that each signal has a broader range of frequencies than the previous ones. The inversion process then begins with using the first signal. Upon convergence, the profile is used as a starting point with the second signal, and the process is repeated for all signals.

### 3.3.3. Biased search direction for $\lambda$

Simultaneous inversion for both $\lambda$ and $\mu$ is remarkably challenging [11]. As we demonstrate in Section 4.2, the objective functional (7) is more sensitive to $\mu$, than to $\lambda$. Consequently, as the inversion evolves, the $\mu$-profile converges faster than that of $\lambda$. In [22], a biasing scheme was proposed to accelerate the convergence of the $\lambda$-profile, such that, at the early stages of inversion, the search direction for $\lambda$ is biased according to that of $\mu$.

The main idea is that due to physical considerations, the $\lambda$-profile should be, more or less, similar to the $\mu$-profile. Hence, during the early inversion iterations, the search direction for $\lambda$ is biased according to:

---

[9]Since we normalize the regularization-part and misfit-part of the gradient in (23b), $\wp = 0.5$ means the gradient is weighted twice by the misfit than the regularization.

$$\boldsymbol{s}_k^\lambda \leftarrow \|\boldsymbol{s}_k^\lambda\| \left( W \frac{\boldsymbol{s}_k^\mu}{\|\boldsymbol{s}_k^\mu\|} + (1-W) \frac{\boldsymbol{s}_k^\lambda}{\|\boldsymbol{s}_k^\lambda\|} \right), \tag{25}$$

where $W$ is a weight that imposes the biasing amount. We assign full weight ($W = 1$) on $\mu$ at the first inversion iteration, and reduce it linearly down to zero as iterates evolve (say at $k = 50$). After that, we let $\lambda$ evolve on its own, according to the original, unbiased search direction.

## 4. Numerical experiments

We present numerical experiments[10] with increasing complexity to test the proposed inversion scheme. In the first example, we verify the accuracy of the gradients, computed by using Algorithm 2. Next, we focus on material profile reconstruction for heterogeneous hosts, using synthetic data at measurement locations. Specifically, we consider: (a) a medium with smoothly varying material properties along depth, to study various aspects of the inversion scheme; (b) a horizontally layered profile with sharp interfaces; (c) a horizontally layered profile with an ellipsoidal inclusion, using highly noisy data; and (d) a layered profile with three inclusions in an attempt to implicate arbitrary heterogeneity. Throughout, we use Gaussian pulses to probe the considered domains:

$$f(t) = e^{-(\frac{t-\bar{\mu}}{\bar{\sigma}})^2},$$

where the parameters that characterize the load are given in Table 1; $\bar{\mu}$ is the mean, $\bar{\sigma}$ is the deviation, $f_{max}$ is the maximal frequency content of the pulse, $t_{end}$ is the active duration of the Gaussian pulse, and the load has an amplitude of 1 kPa. The time history of the loads and their corresponding Fourier spectrum are shown in Figure 4.

Table 1: Characterization of Gaussian pulses.

| Name | $f_{max}$ | $\bar{\mu}$ | $\bar{\sigma}$ | $t_{end}$ |
|------|-----------|-------------|----------------|-----------|
| $p_{20}$ | 20 | 0.11 | 0.0014 | 0.20 |
| $p_{30}$ | 30 | 0.08 | 0.0007 | 0.15 |
| $p_{40}$ | 40 | 0.06 | 0.0004 | 0.12 |

*4.1. Numerical verification of the material gradients*

Accurate computation of the discrete gradients is crucial for the robustness of Algorithm 2. The gradients of the objective functional with respect to the control parameters can be computed either by the optimize-then-discretize, or, the discretize-then-optimize approach [30]. While the discretize-then-optimize method yields the exact discrete gradients of the discrete objective functional [16], this is not always the case with the optimize-then-discretize scheme [43].

---

[10] We developed a code in Fortran, using PETSc [3] to facilitate parallel implementation.

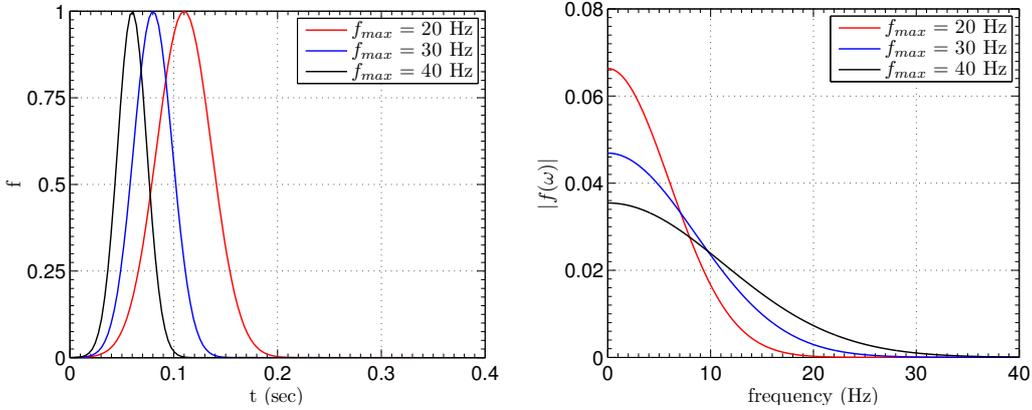

Figure 4: Time history of the Gaussian pulses and their Fourier spectrum.

In this part, through a numerical experiment, we demonstrate that the discrete gradients that we compute via the optimize-then-discretize technique, are accurate, and equal to the discrete gradients of the discrete objective functional. To this end, we compare directional finite differences of the discrete objective functional, with directional gradients obtained from (18). We start by defining the finite difference directional derivatives:

$$d_h^{\mathrm{fd}}(\boldsymbol{\lambda},\boldsymbol{\mu})(\tilde{\boldsymbol{\lambda}}) := \frac{\mathbf{J}(\boldsymbol{\lambda}+h\tilde{\boldsymbol{\lambda}},\boldsymbol{\mu}) - \mathbf{J}(\boldsymbol{\lambda},\boldsymbol{\mu})}{h}, \tag{26a}$$

$$d_h^{\mathrm{fd}}(\boldsymbol{\lambda},\boldsymbol{\mu})(\tilde{\boldsymbol{\mu}}) := \frac{\mathbf{J}(\boldsymbol{\lambda},\boldsymbol{\mu}+h\tilde{\boldsymbol{\mu}}) - \mathbf{J}(\boldsymbol{\lambda},\boldsymbol{\mu})}{h}, \tag{26b}$$

where $\tilde{\boldsymbol{\lambda}}$ and $\tilde{\boldsymbol{\mu}}$ is the discrete direction vector for $\lambda$ and $\mu$, respectively. The directional derivatives obtained via the control problems (18) are:

$$d^{\mathrm{co}}(\boldsymbol{\lambda},\boldsymbol{\mu})(\tilde{\boldsymbol{\lambda}}) = \tilde{\boldsymbol{\lambda}}^T \, \tilde{\mathbf{M}} \, \mathbf{g}^\lambda, \tag{27a}$$

$$d^{\mathrm{co}}(\boldsymbol{\lambda},\boldsymbol{\mu})(\tilde{\boldsymbol{\mu}}) = \tilde{\boldsymbol{\mu}}^T \, \tilde{\mathbf{M}} \, \mathbf{g}^\mu. \tag{27b}$$

Next, we verify that (26) and (27) produce identical values for several choices that we make for $\tilde{\boldsymbol{\lambda}}$ and $\tilde{\boldsymbol{\mu}}$, by considering a test problem displayed in Figure 5, and detailed below, with regularization factors $R_\lambda = R_\mu = 0$[11]. We consider perturbations $\tilde{\boldsymbol{\lambda}}$ or $\tilde{\boldsymbol{\mu}}$: the unit

---

[11]Zero values are considered since convergence difficulties that may arise stem from the misfit part of the

vector is zero everywhere except at the component corresponding to coordinate $(x, y, z)$ where the directional derivatives are being computed. The derivatives $d^{co}$ and $d_h^{fd}$ with respect to either $\lambda$ or $\mu$, for points with coordinates $(x, y, z)$, are presented in Table 2. Digits where $d_h^{fd}$ coincides with $d^{co}$ are shown in bold. Since pointwise perturbations result in small changes in the objective functional, numerical roundoff influences the accuracy of the finite difference directional derivatives, as it has also been reported in [43].

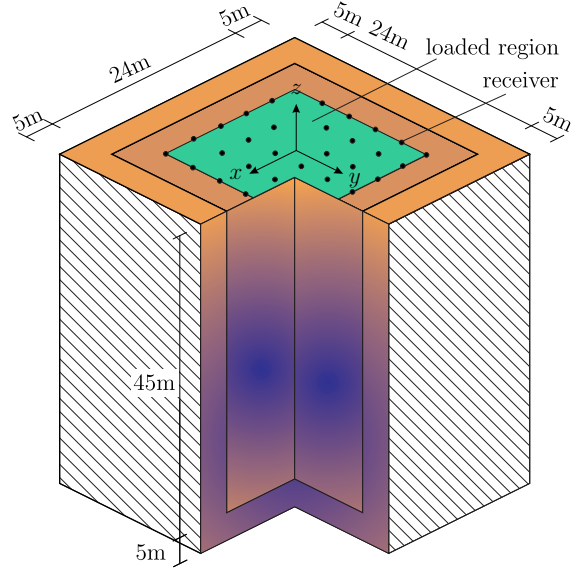

Figure 5: Problem configuration for the verification of the gradients.

Table 2: Comparison of the directional derivatives.

| Case | $f_{max}$ | (x,y,z) | Pert. field | $d^{co}$ | $d_h^{fd}$ $h = 10^{-3}$ | $h = 10^{-4}$ | $h = 10^{-5}$ |
|---|---|---|---|---|---|---|---|
| 1 | 20 Hz | (1,1,0) | $\lambda$ | $-3.03500$e-9 | $-\mathbf{3.03}416$e-9 | $-\mathbf{3.03}496$e-9 | $-\mathbf{3.03501}$e-9 |
| 2 | 20 Hz | (1,1,0) | $\mu$ | $-2.78908$e-9 | $-\mathbf{2.78}875$e-9 | $-\mathbf{2.789}17$e-9 | $-\mathbf{2.789}21$e-9 |
| 3 | 20 Hz | (1,1,-40) | $\lambda$ | $-5.14848$e-11 | $-\mathbf{5.14}711$e-11 | $-\mathbf{5.14}647$e-11 | $-\mathbf{5.14}996$e-11 |
| 4 | 20 Hz | (1,1,-40) | $\mu$ | $+4.97666$e-10 | $+\mathbf{4.97}411$e-10 | $+\mathbf{4.97}512$e-10 | $+\mathbf{4.97}366$e-10 |
| 5 | 40 Hz | (1,1,0) | $\lambda$ | $-1.07645$e-9 | $-\mathbf{1.076}23$e-9 | $-\mathbf{1.0765}2$e-9 | $-\mathbf{1.0765}6$e-9 |
| 6 | 40 Hz | (1,1,0) | $\mu$ | $-1.56155$e-9 | $-\mathbf{1.5615}3$e-9 | $-\mathbf{1.561}78$e-9 | $-\mathbf{1.561}80$e-9 |

We remark that the agreement between the two derivatives is remarkable, both for cases 1-4, where the wavefield is well-resolved, and for cases 5 and 6, where only 10 points

---

objective functional, and not from the regularization part. Nevertheless, we have also successfully verified the accuracy of the regularization component of the gradients.

per wavelength are used for spatial discretization.

The considered test problem is a heterogeneous half-space with a smoothly varying material profile along depth, given in (28), and mass density $\rho = 2000$ kg/m$^3$, which, after truncation, is reduced to a cubic computational domain of length and width 24 m × 24 m, and 45 m depth. A 5 m-thick PML is placed at the truncation boundaries, as shown in Figure 5. The material properties at the interfaces $\Gamma^\text{I}$ are extended into the PML. The interior and PML domains are discretized by quadratic hexahedral spectral elements (i.e., 27-noded bricks, and quadratic-quadratic pairs of approximation for displacement and stress components in the PML, and, also, quadratic approximation for material properties) of size 1 m, and $\Delta t = 9 \times 10^{-4}$ s. Throughout, for the PML parameters, we choose $\alpha_o = 5$, $\beta_o = 400$ s$^{-1}$, and a quadratic profile for the attenuation functions, i.e., $m = 2$. See [13] for notation and other details. To probe the medium, we consider vertical stress loads with Gaussian pulse temporal signatures (see Table 1), applied on the surface of the domain over a region $(-11 \text{ m} \leq x, y \leq 11 \text{ m})$, whereas receivers that collect displacement response $\mathbf{u}_m(\mathbf{x}, t)$ are also located in the same region, at every grid point. To obtain synthetic data at the receiver locations, we use a model with identical characteristics and dimensions as detailed above, but, with a refined discretization; i.e., element size of 0.5 m, and $\Delta t = 4.5 \times 10^{-4}$ s. The data was then mapped onto the coarser mesh discussed earlier. The total duration of the simulation is $T = 0.5$ s.

### 4.2. Smoothly varying heterogeneous medium

We consider a heterogeneous half-space with a smoothly varying material profile along depth, given by:

$$\lambda(z) = \mu(z) = 80 + 0.45 \, |z| + 35 \, \exp\left(-\frac{(|z| - 22.5)^2}{150}\right) \text{ (MPa)}, \tag{28}$$

and mass density $\rho = 2000$ kg/m$^3$, which, after truncation, is reduced to a cubic computational domain of length and width 40 m × 40 m, and 45 m depth. A 6.25 m-thick PML is placed at the truncation boundaries, as illustrated in Figure 6. The target profiles are shown in Fig, 7. The material properties at the interfaces $\Gamma^\text{I}$ are extended into the PML. The interior and PML domains are discretized by quadratic hexahedral spectral elements (i.e., 27-noded bricks, and quadratic-quadratic pairs of approximation for displacement and stress components in the PML, and, also, quadratic approximation for material properties) of size 1.25 m, and $\Delta t = 10^{-3}$ s. This leads to $3,578,136$ state unknowns, and $616,850$ material parameters. To probe the medium, we consider vertical stress loads with Gaussian pulse temporal signatures (see Table 1), applied on the surface of the domain over a region $(-17.5 \text{ m} \leq x, y \leq 17.5 \text{ m})$, whereas receivers that collect displacement response $\mathbf{u}_m(\mathbf{x}, t)$ are placed at every grid point, in the same region.

Before attempting simultaneous inversion for the two Lamé parameters, we perform single parameter inversion for a) $\mu$ only, assuming $\lambda$ is *a priori* known and fixing it to the target profile; and b) $\lambda$ only, assuming distribution of $\mu$ is known.

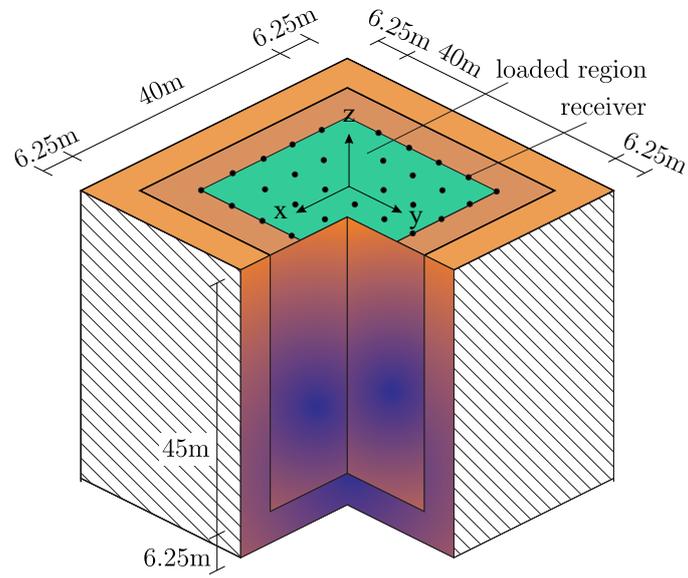

Figure 6: Problem configuration: material profile reconstruction of a smoothly varying medium.

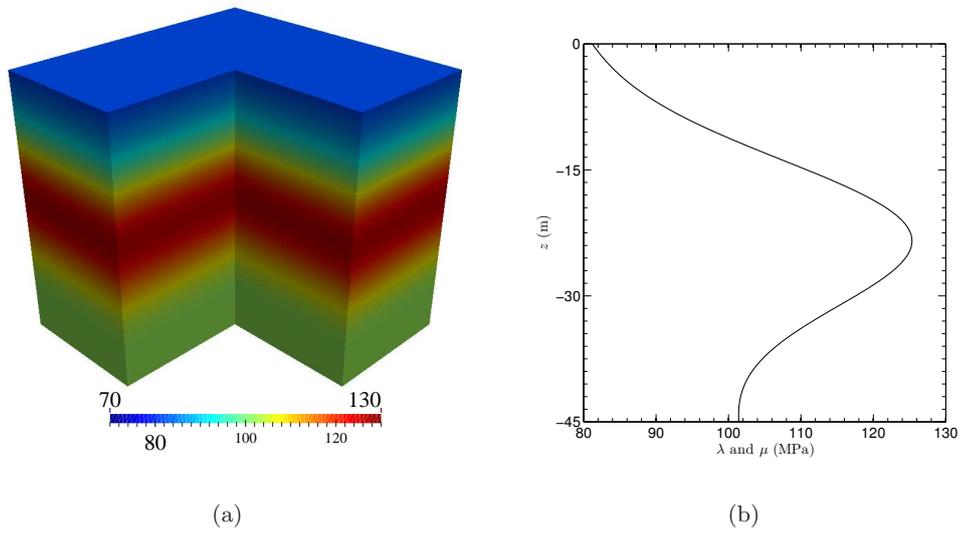

(a)

(b)

Figure 7: Smoothly varying medium: (a) target $\lambda$ and $\mu$ (MPa); and (b) profile at $(x, y) = (0, 0)$.

*4.2.1. Single parameter inversion*

First, we assume $\lambda$ is *a priori* known, and fix it to the target profile. We start inverting for $\mu$, with a homogeneous initial guess of 80 MPa, exploiting Tikhonov regularization for taming ill-posedness and solution multiplicity. We use the Gaussian pulse $p_{20}$ with maximal frequency content $f_{max} = 20$ Hz (see Table 1) for 50 iterations, and, then, switch to $p_{30}$ with $f_{max} = 30$ Hz. After 156 iterations, $\mu$ converges to the target profile, as shown in Figure 8(b). We compare the inverted cross-sectional profiles of $\mu$ with the target profile at three different locations, as shown in Figure 8(c), 8(d), and 8(e). The agreement between the two profiles is excellent. Reduction of the misfit functional with respect to inversion iterations is shown in Figure 10(a), which is almost 7 orders of magnitude.

Next, we fix $\mu$ to the target profile, and invert for $\lambda$, starting with a homogeneous initial guess of 80 MPa. We use the Gaussian pulse $p_{20}$ for 160 iterations, $p_{30}$ up to the $300^{th}$ iteration, and then switch to $p_{40}$. After 456 iterations, the optimizer converges to the profile displayed in Figure 9(a). The agreement between the target profile and the inverted profile is remarkable. We compare the two profiles at three different cross-sections shown in Figure 9(c), 9(d), and 9(e): the agreement between the two profiles is excellent. The misfit history is shown in Figure 10(b); the optimizer reduced the misfit almost 6 orders of magnitude.

We remark that the initial value of the misfit in the first experiment is almost 2 orders of magnitude more than that of the second experiment. This indicates that the objective functional is not equally sensitive to both control parameters, as it has also been reported in [22]: the objective functional is more sensitive to $\mu$.

*4.2.2. Simultaneous inversion*

We start with a homogeneous initial guess of 80 MPa for both $\lambda$ and $\mu$ and attempt simultaneous inversion. The target profiles are shown in Figure 7, and the inverted profiles are displayed in Figure 11(a) and 11(b). We also compare the cross-sectional values of the target and inverted profiles at three different locations, shown in Figure 12. Although the inverted $\mu$ profile agrees reasonably well with the target profile, inversion for $\lambda$ is not satisfactory, and the inverted profile departs from the target as depth increases.

Due to the unsuccessful inversion of the $\lambda$ profile in the case of simultaneous inversion, in the next experiment, we bias the search direction of $\lambda$ based on that of $\mu$, at the very early stages of inversion, according to the procedure detailed in Section 3.3.3. This leads to the successful reconstruction of the two profiles, as is shown in Figure 13(a) and 13(b). In Figure 14, we compare the cross-sectional values of the target and the inverted profiles. The agreement of the inverted $\mu$ profile with the target is remarkable. Moreover, the inverted $\lambda$ profile agrees reasonably well with the target, with some discrepancies in depth. The misfit history is shown in Figure 15(b), where the kink in the misfit curve at the $50^{th}$ iteration corresponds to the termination point of the biasing scheme.

We remark that in practical applications, one is more interested in the shear wave velocity ($c_s$) and compression wave velocity ($c_p$) profiles. Once the Lamé parameters have

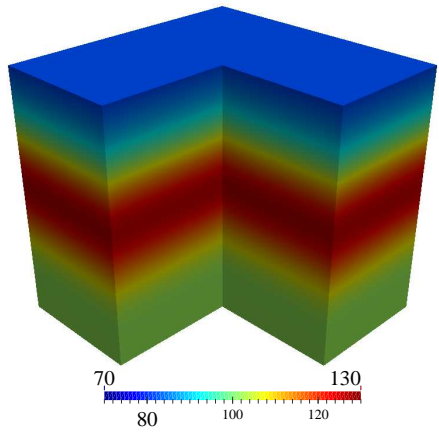
(a) $\lambda$ (*a priori* known)

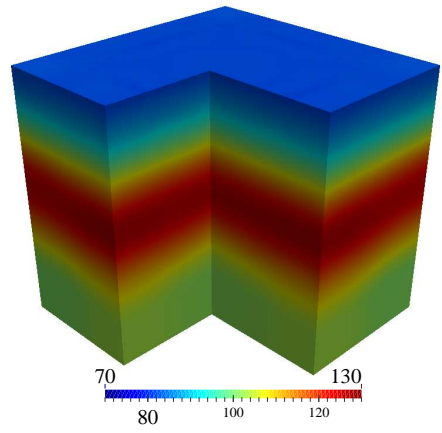
(b) $\mu$ (inverted)

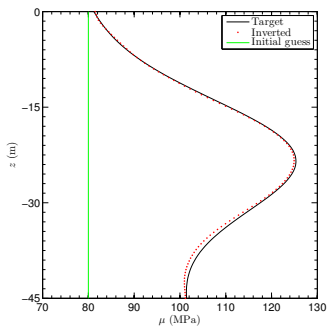
(c) $\mu$ at $(x, y) = (0, 0)$

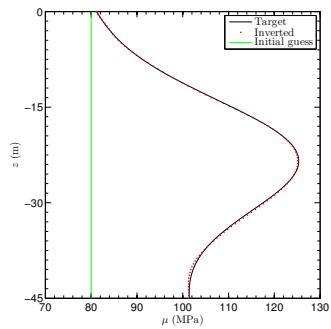
(d) $\mu$ at $(x, y) = (10, 10)$

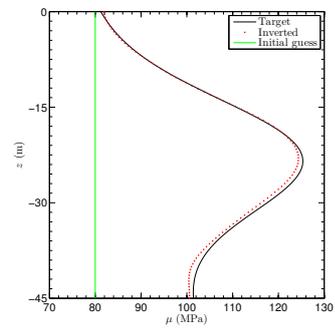
(e) $\mu$ at $(x, y) = (20, 20)$

Figure 8: Single-parameter inversion ($\mu$ only) for a smoothly varying medium.

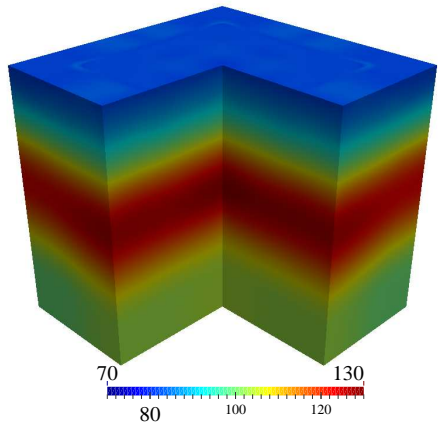

(a) $\lambda$ (inverted)

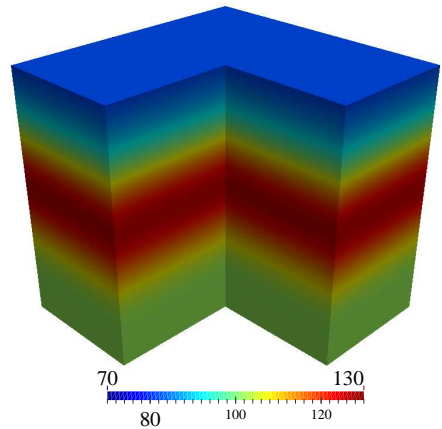

(b) $\mu$ (*a priori* known)

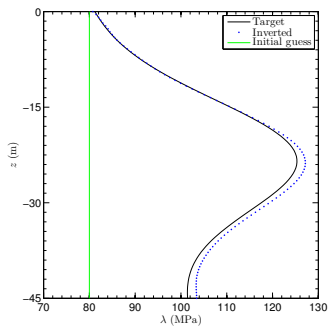

(c) $\lambda$ at $(x, y) = (0, 0)$

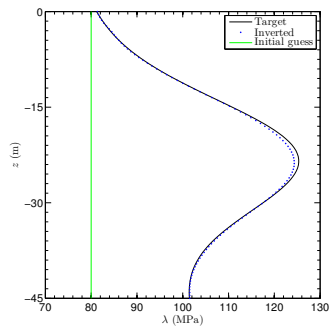

(d) $\lambda$ at $(x, y) = (10, 10)$

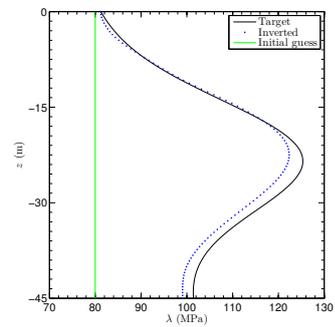

(e) $\lambda$ at $(x, y) = (20, 20)$

Figure 9: Single-parameter inversion ($\lambda$ only) for a smoothly varying medium.

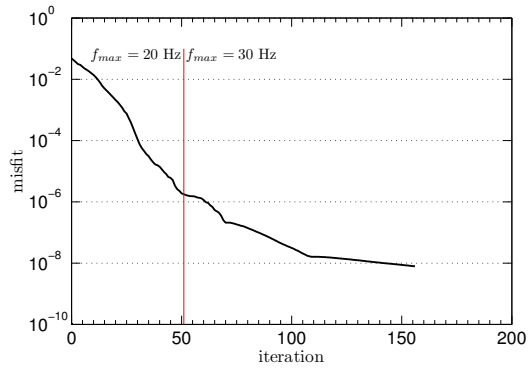
(a) reconstruction for $\mu$ only

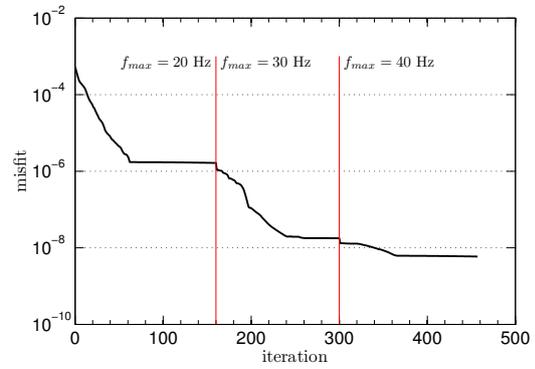
(b) reconstruction for $\lambda$ only

Figure 10: Variation of the misfit functional with respect to inversion iterations (single parameter inversion).

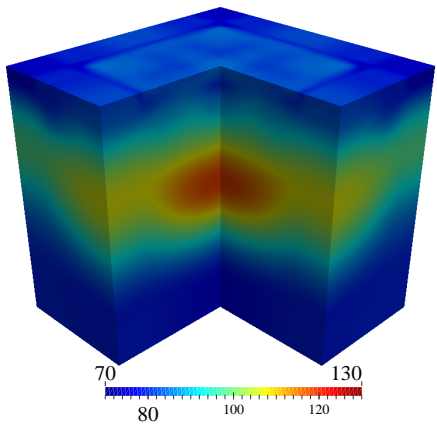
(a) $\lambda$ (inverted)

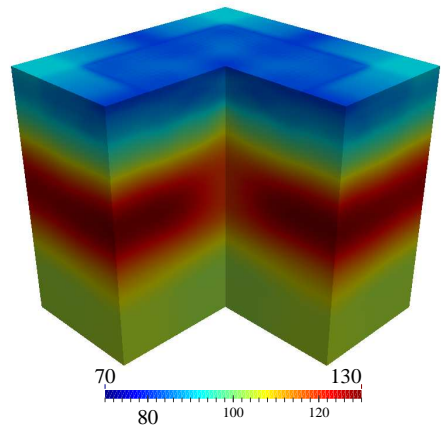
(b) $\mu$ (inverted)

Figure 11: Simultaneous inversion for $\lambda$ and $\mu$ using unbiased search directions.

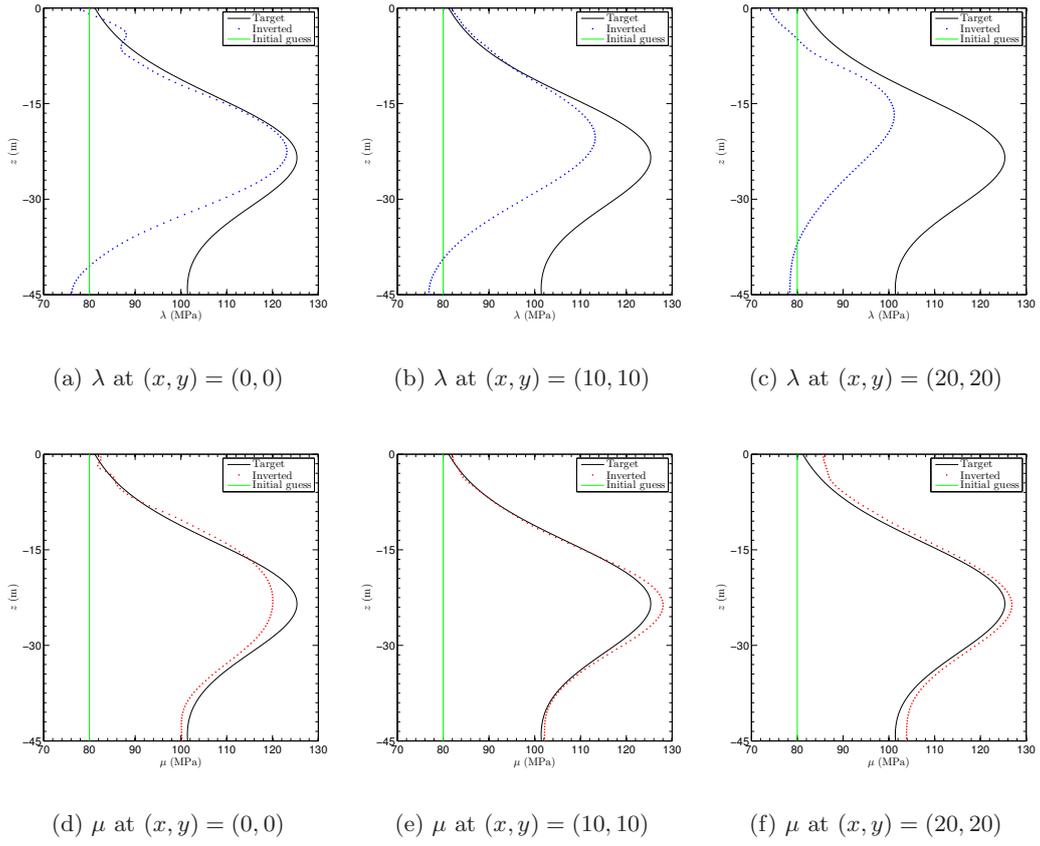

(a) $\lambda$ at $(x,y) = (0,0)$  (b) $\lambda$ at $(x,y) = (10,10)$  (c) $\lambda$ at $(x,y) = (20,20)$

(d) $\mu$ at $(x,y) = (0,0)$  (e) $\mu$ at $(x,y) = (10,10)$  (f) $\mu$ at $(x,y) = (20,20)$

Figure 12: Cross-sectional profiles of $\lambda$ and $\mu$ using unbiased search directions.

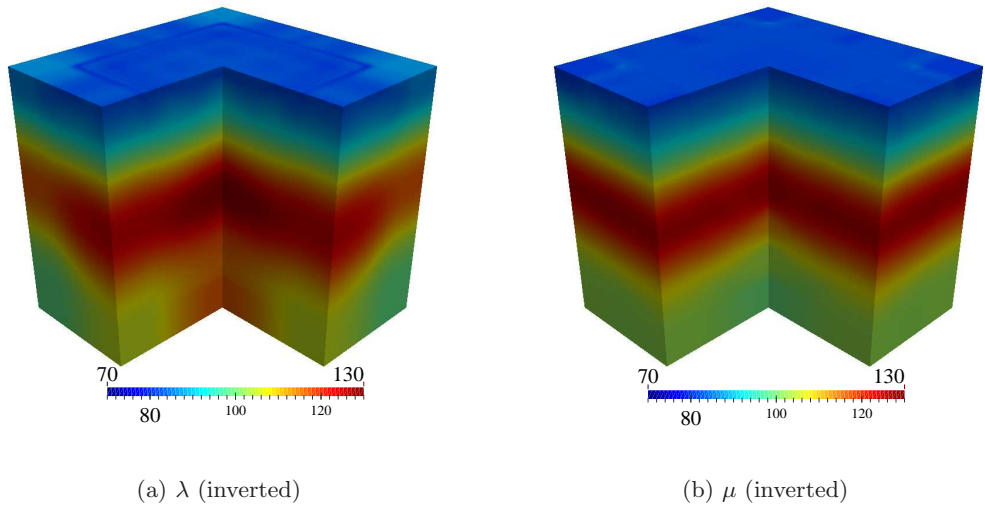

(a) $\lambda$ (inverted)  (b) $\mu$ (inverted)

Figure 13: Simultaneous inversion for $\lambda$ and $\mu$ using biased search directions.

been determined, the wave velocities can be readily computed via:

$$c_s = \sqrt{\frac{\mu}{\rho}}, \qquad c_p = \sqrt{\frac{\lambda + 2\mu}{\rho}}. \tag{29}$$

In Figure 16, we compare the compression wave velocities at three different cross-sectional locations, where the agreement between the reconstructed $c_p$ profile and the target is remarkable. The shear wave velocity does not depend on $\lambda$, and, therefore, its quality is similar to that of the $\mu$ profile.

4.3. Layered medium

We consider a 40 m × 40 m × 45 m layered medium, where a 6.25 m-thick PML is placed at its truncation boundaries. The properties of the medium are:

$$\lambda(z) = \mu(z) = \begin{cases} 80 \text{ MPa}, & \text{for } -12 \text{ m} \leq z \leq 0 \text{ m}, \\ 101.25 \text{ MPa}, & \text{for } -27 \text{ m} \leq z < -12 \text{ m}, \\ 125 \text{ MPa}, & \text{for } -50 \text{ m} \leq z < -27 \text{ m}, \end{cases} \tag{30}$$

and are shown in Figure 17, with mass density $\rho = 2000$ kg/m$^3$. The material properties at the interfaces $\Gamma^{\text{I}}$ are extended into the PML buffer. The interior and PML domains are discretized by quadratic hexahedral spectral elements (i.e., 27-noded bricks, and quadratic-quadratic pairs of approximation for displacement and stress components in the PML, and,

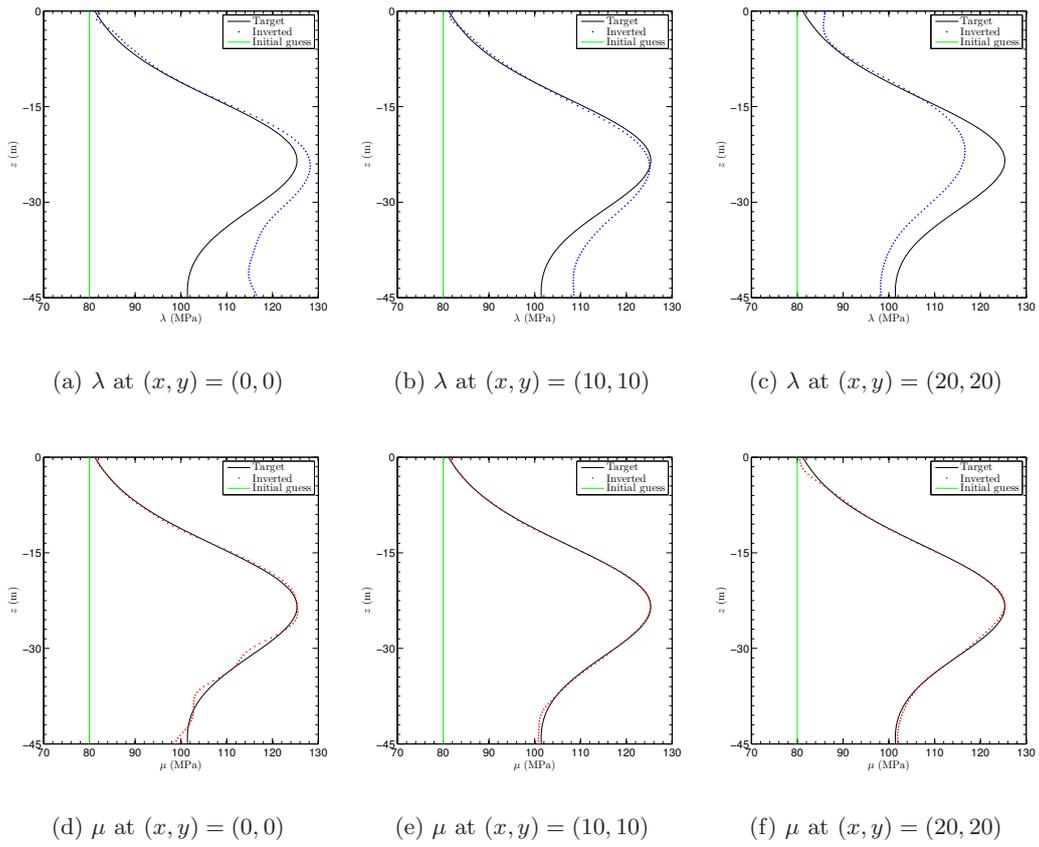

(a) $\lambda$ at $(x,y) = (0,0)$  (b) $\lambda$ at $(x,y) = (10,10)$  (c) $\lambda$ at $(x,y) = (20,20)$

(d) $\mu$ at $(x,y) = (0,0)$  (e) $\mu$ at $(x,y) = (10,10)$  (f) $\mu$ at $(x,y) = (20,20)$

Figure 14: Cross-sectional profiles of $\lambda$ and $\mu$ using biased search directions.

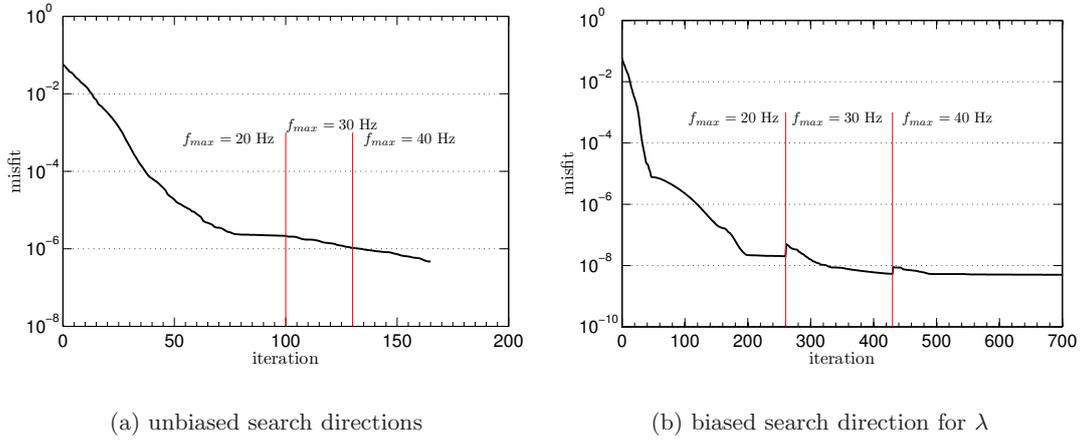

(a) unbiased search directions

(b) biased search direction for $\lambda$

Figure 15: Variation of the misfit functional with respect to inversion iterations (simultaneous inversion).

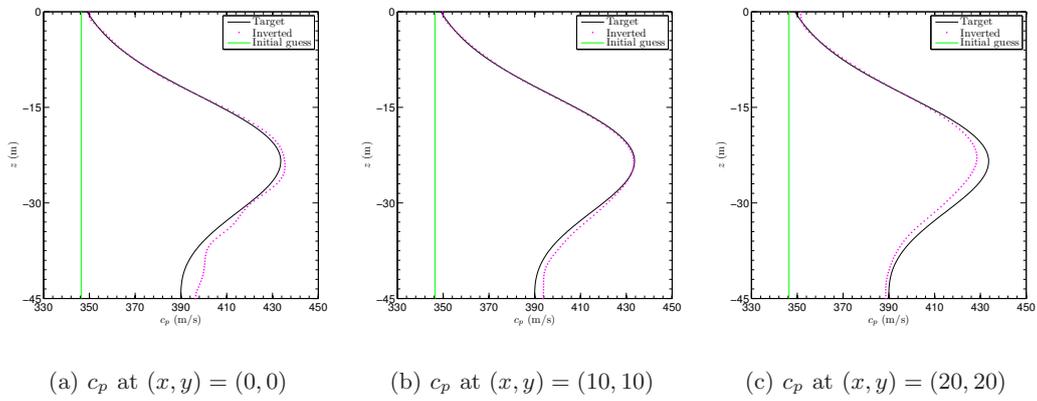

(a) $c_p$ at $(x,y) = (0,0)$

(b) $c_p$ at $(x,y) = (10,10)$

(c) $c_p$ at $(x,y) = (20,20)$

Figure 16: Cross-sectional profiles of $c_p$ using biased search directions.

also, quadratic approximation for material properties) of size 1.25 m, and $\Delta t = 10^{-3}$ s. For probing the medium, we use vertical stress loads with Gaussian pulse temporal signatures, applied on the surface of the domain over a region $(-17.5\text{m} \leq x, y \leq 17.5\text{m})$, whereas receivers that collect displacement response $\mathbf{u}_m(\mathbf{x}, t)$ are also located in the same region, at every grid point.

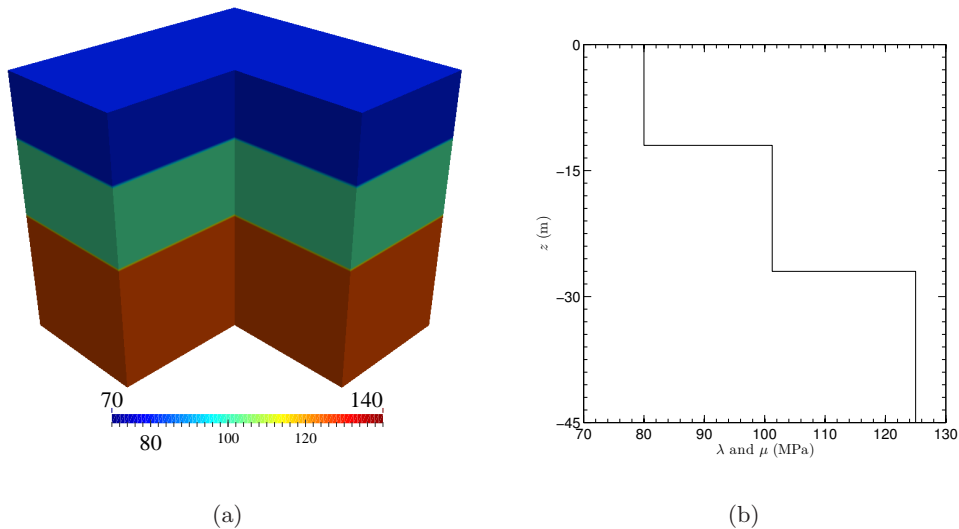

Figure 17: Layered medium: (a) target $\lambda$ and $\mu$ (MPa); and (b) profile at $(x, y) = (0, 0)$.

We start the inversion process with a homogeneous initial guess of 80 MPa for the Lamé parameters, and attempt simultaneous inversion for both $\lambda$ and $\mu$, using the biasing scheme outlined in Section 3.3.3. We use the Total Variation regularization scheme, with $\epsilon = 0.01$, to capture the sharp interfaces of the target profiles. We use the Gaussian pulse $p_{20}$, with $f_{max} = 20$ Hz, and final simulation time[12] $T = 0.45$ s, for 310 iterations. The resulting profiles are shown in Figure 18(a) and 18(b), where the layering of the medium is clearly visible in the inverted profiles. To improve the quality of the inverted profiles, we use them as an initial guess with the Gaussian pulse $p_{30}$, and final simulation time of $T = 0.4$ s, for up to the $860^{\text{th}}$ iteration, and, then, switch to $p_{40}$, with final simulation time of $T = 0.4$ s. After 1112 iterations, the optimizer converges to the profiles displayed in

---

[12]The final simulation time is chosen such that they are long enough to probe the medium effectively, but not too long to increse the computational cost, unnecessarily. To this end, we run a forward simulation by using the target profile, and monitor when the total energy of the system (see [13]) die out. We then use this time duration for the inversion process. Reference [20] provides guidelines for choosing the source duration when inversion is performed by using field data.

Figure 18(c) and 18(d). There is excellent agreement between the inverted $\mu$ profile and the target profile. The inverted $\lambda$ profile is also in good agreement with the target profile: the two top layers have been reconstructed quite well, whereas the bottom layer is slightly "stiffer" in it's middle zone. We compare the inverted profiles with the targets at three different cross-sections, shown in Figure 19. Due to the TV regularization, sharp interfaces have been captured quite successfully. In Figure 20, we compare the $c_p$ profile with the target, at the same cross-sections; the agreement between the two profiles is impressive. Figure 21 shows the misfit history: the optimizer reduced the misfit almost 7 orders of magnitude.

### 4.4. Layered medium with inclusion

We consider a layered medium with an inclusion. The problem consists of a 40 m × 40 m × 45 m layered medium with an ellipsoidal inclusion, where a 6.25 m-thick PML is placed at its truncation boundaries. The material profiles are given by:

$$\lambda(z) = \mu(z) = \begin{cases} 80 \text{ MPa}, & \text{for } -12 \text{ m} \leq z \leq 0 \text{ m}, \\ 101.25 \text{ MPa}, & \text{for } -27 \text{ m} \leq z < -12 \text{ m}, \\ 125 \text{ MPa}, & \text{for } -50 \text{ m} \leq z < -27 \text{ m}, \\ 156.8 \text{ MPa}, & \text{for ellipsoidal inclusion}, \end{cases} \quad (31)$$

and are shown in Figure 22, with constant mass density $\rho = 2000$ kg/m$^3$. The ellipsoidal inclusion occupies the region $(\frac{x-7.5}{7.5})^2 + (\frac{y}{5})^2 + (\frac{z+12}{5.5})^2 \leq 1$. The material properties at the interfaces $\Gamma^I$ are extended into the PML buffer. The interior and PML domains are discretized by quadratic hexahedral spectral elements (i.e., 27-noded bricks, and quadratic-quadratic pairs of approximation for displacement and stress components in the PML, and, also, quadratic approximation for material properties) of size 1.25 m, and $\Delta t = 10^{-3}$ s. To illuminate the domain, we use vertical stress loads with Gaussian pulse temporal signatures, applied on the surface of the medium over a region ($-17.5$m $\leq x, y \leq 17.5$m), whereas the receivers are also placed at every grid point in the same region.

We use the Total Variation regularization scheme to alleviate ill-posedness and solution multiplicity, with $\epsilon = 0.01$. Similar to the previous examples, we use a source-frequency continuation scheme, starting with the Gaussian pulse $p_{20}$ with maximal frequency content of 20 Hz for $T = 0.45$ s, and, when updates in the material profiles become practically insignificant, we switch to the next load in Table 1, which contains a broader range of frequencies, and, therefore, is able to image finer features. Figure 23(a) and 23(b) show the material profiles after 410 iterations, which adequately capture the layering of the domain as well as the ellipsoidal inclusion. To improve the quality of the reconstructed profiles, we use them as an initial guess with the Gaussian pulse $p_{30}$, and final simulation time of $T = 0.4$ s, for up to the 730$^{\text{th}}$ iteration, and, then, switch to $p_{40}$, with final simulation time of $T = 0.4$ s. Figure 23(c) and 23(d) show the inverted profiles after 1160 iterations. The sharp interfaces between the three layers and around the ellipsoidal inclusion are very well

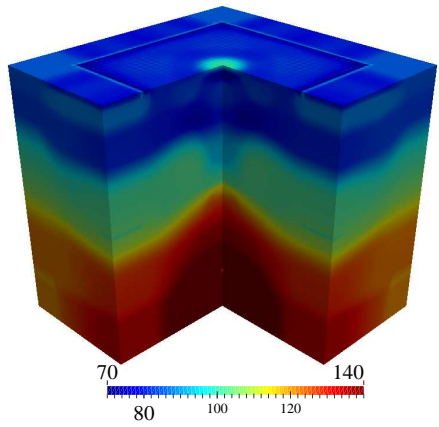

(a) $\lambda$ ($f_{max} = 10$ Hz)

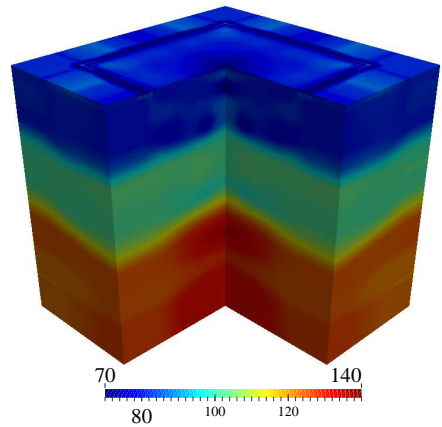

(b) $\mu$ ($f_{max} = 10$ Hz)

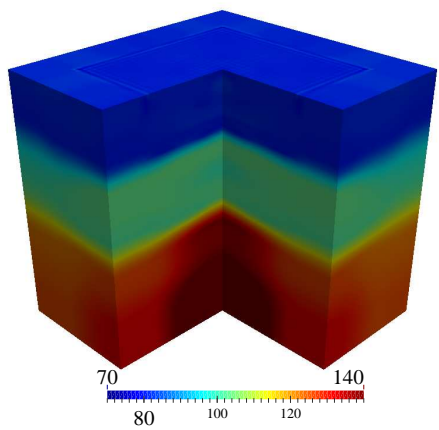

(c) $\lambda$ ($f_{max} = 40$ Hz)

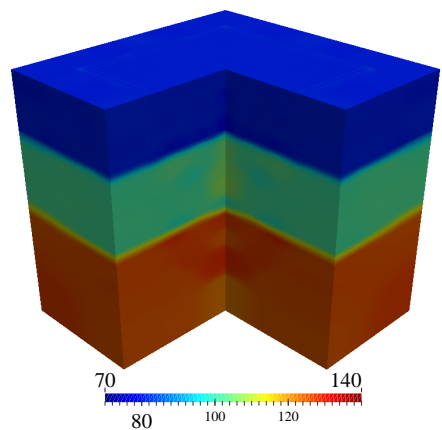

(d) $\mu$ ($f_{max} = 40$ Hz)

Figure 18: Simultaneous inversion for $\lambda$ and $\mu$ (layered medium).

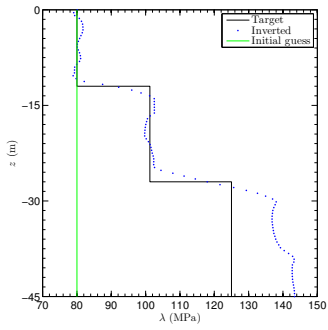 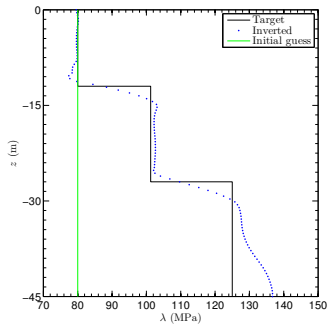 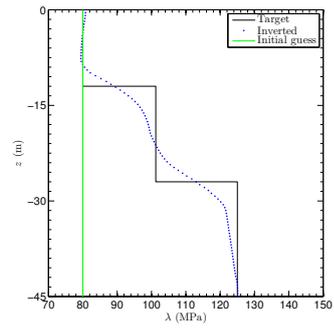

(a) $\lambda$ at $(x,y)=(0,0)$     (b) $\lambda$ at $(x,y)=(10,10)$     (c) $\lambda$ at $(x,y)=(20,20)$

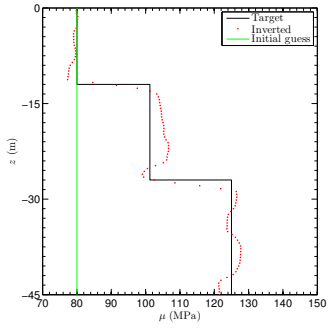 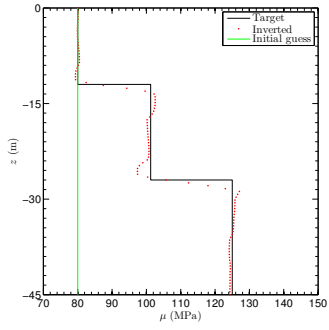 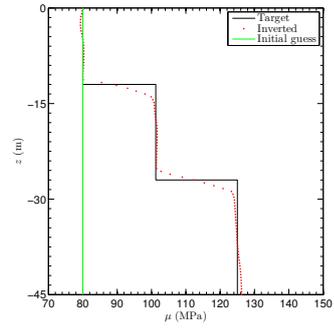

(d) $\mu$ at $(x,y)=(0,0)$     (e) $\mu$ at $(x,y)=(10,10)$     (f) $\mu$ at $(x,y)=(20,20)$

Figure 19: Cross-sectional profiles of $\lambda$ and $\mu$ (layered medium).

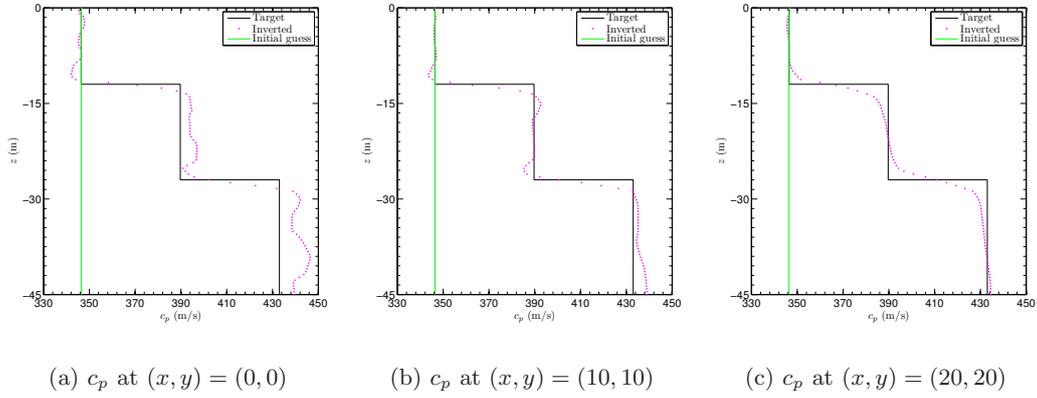

(a) $c_p$ at $(x,y) = (0,0)$  (b) $c_p$ at $(x,y) = (10,10)$  (c) $c_p$ at $(x,y) = (20,20)$

Figure 20: Cross-sectional profiles of $c_p$ (layered medium).

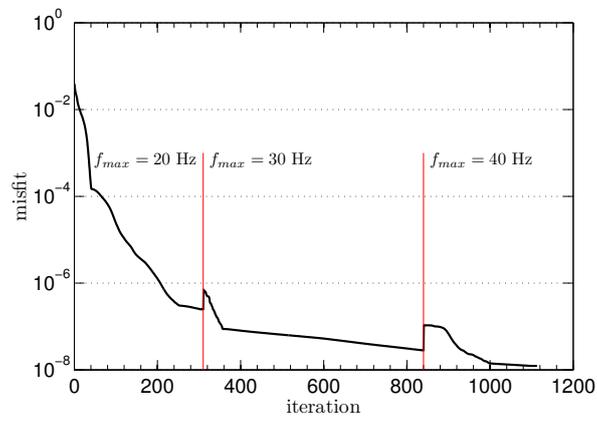

Figure 21: Variation of the misfit functional with respect to inversion iterations (layered medium).

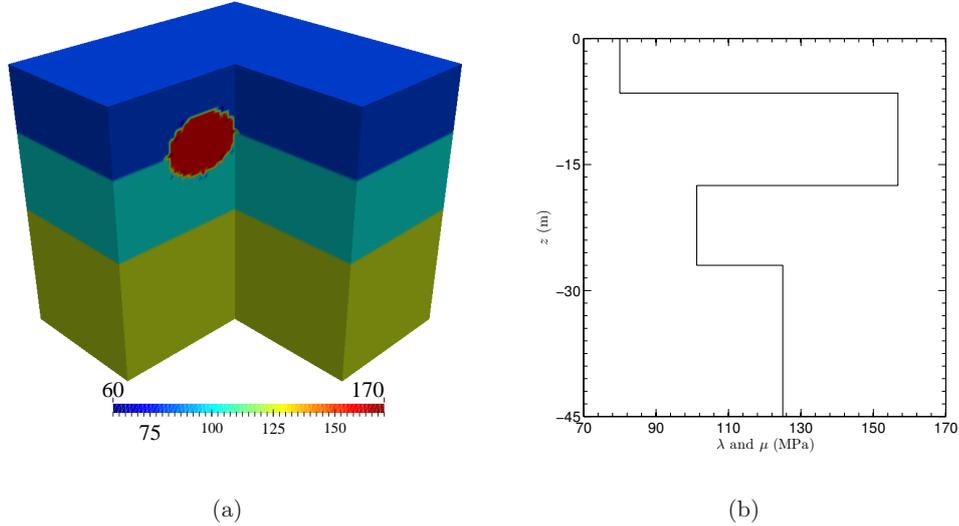

Figure 22: Layered medium with inclusion: (a) target $\lambda$ and $\mu$; and (b) profile at $(x,y) = (7.5, 0)$.

captured for the $\mu$ profile. The $\lambda$ profile agrees reasonably well with the target, showing some "stiff" features at the center of the bottom layer, similar to the previous example.

Figures 24 and 25 compare the inverted profiles with the target profiles at three different cross-sectional lines of the domain, indicating successful imaging of both the layering and the inclusion. Variation of the misfit functional with respect to the inversion iterations is shown in Figure 25, where, again, a kink at the $50^{\text{th}}$ iteration of the misfit curve, corresponds to the termination point of the biasing scheme.

Encouraged by the successful performance of the proposed inversion algorithm with noise-free data, next, we consider adding different levels of Gaussian noise to the measured synthetic response at the receiver locations, and investigate its effect on the inversion. Figures 27(a)-27(d) show the measured displacement response of the system at $(x,y,z) = (3.125, -13.75, 0)$ m, subjected to the $p_{20}$ pulse, contaminated with 1%, 5%, 10%, and 20% Gaussian noise, respectively. Using the source-frequency continuation scheme, the optimizer converges after 811 and 751 iterations, respectively, for cases corresponding to the 1% and 5% Gaussian noise levels. The inverted profiles are shown in Figure 28. The reconstruction is successful, with minor discrepancies on the top surface. Next, we increase the noise level to 10% and 20%, and attempt inversion; after 770 and 674 iterations, respectively, we converge to the profiles shown in Figure 29. The quality of the inverted profiles decreases as the noise level increases. However, similarly to the previous case, except for a thin layer on the top surface, inversion is successful. In Figure 30, we compare cross-sectional profiles of $\lambda$ and $\mu$ with the target, at different noise levels, at $(x,y) =$

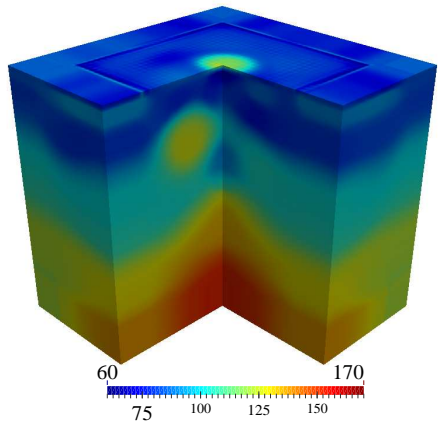
(a) $\lambda$ ($f_{max} = 10$ Hz)

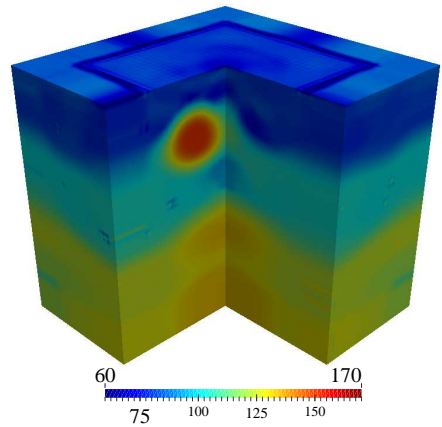
(b) $\mu$ ($f_{max} = 10$ Hz)

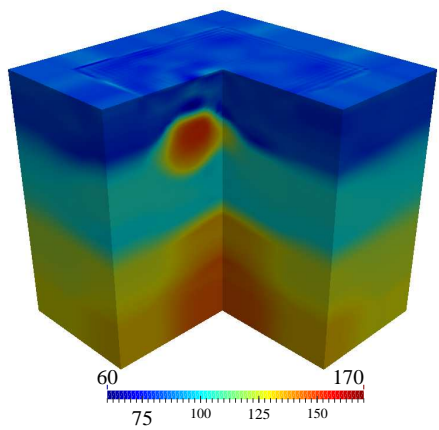
(c) $\lambda$ ($f_{max} = 40$ Hz)

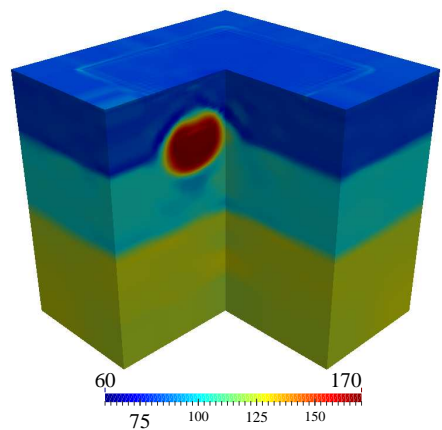
(d) $\mu$ ($f_{max} = 40$ Hz)

Figure 23: Simultaneous inversion for $\lambda$ and $\mu$ (layered medium with inclusion).

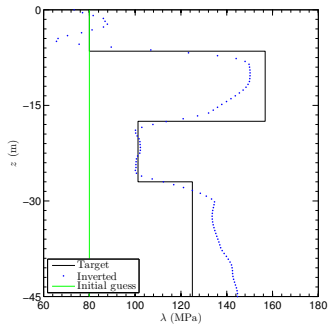
(a) $\lambda$ at $(x,y) = (7.5, 0)$

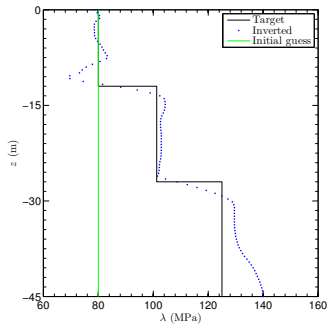
(b) $\lambda$ at $(x,y) = (7.5, 10)$

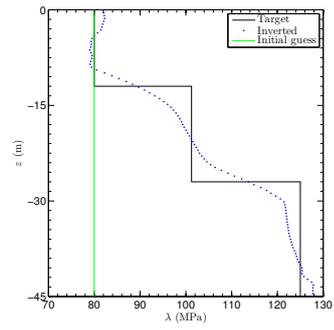
(c) $\lambda$ at $(x,y) = (20, 20)$

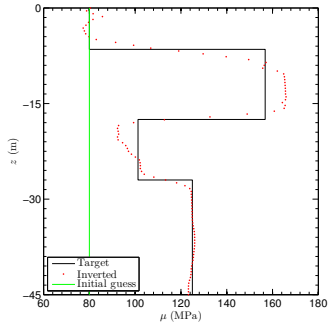
(d) $\mu$ at $(x,y) = (7.5, 0)$

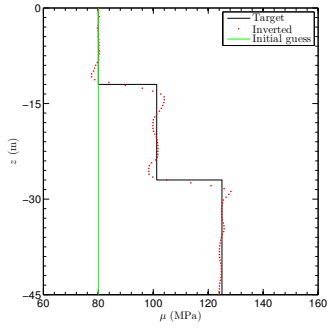
(e) $\mu$ at $(x,y) = (7.5, 10)$

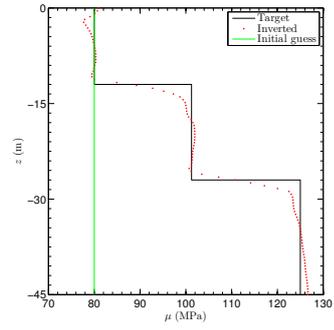
(f) $\mu$ at $(x,y) = (20, 20)$

Figure 24: Cross-sectional profiles of $\lambda$ and $\mu$ (layered medium with inclusion).

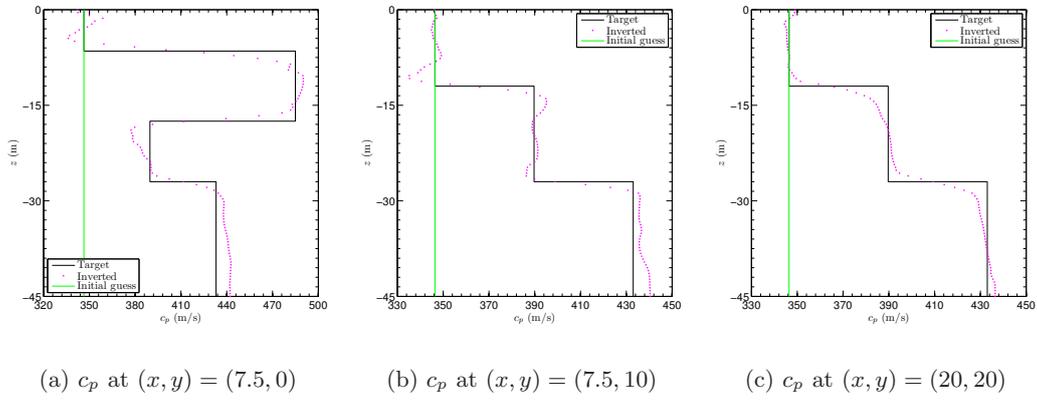

(a) $c_p$ at $(x,y) = (7.5, 0)$  (b) $c_p$ at $(x,y) = (7.5, 10)$  (c) $c_p$ at $(x,y) = (20, 20)$

Figure 25: Cross-sectional profiles of $c_p$ (layered medium with inclusion).

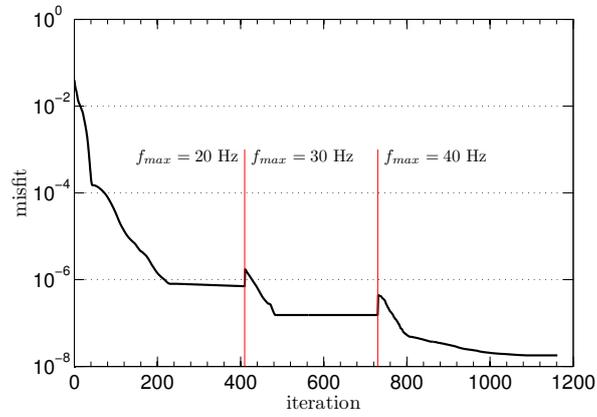

Figure 26: Variation of the misfit functional with respect to inversion iterations (layered medium with inclusion).

$(7.5, 0)$ m, which passes through the center of the ellipsoidal inclusion. Sharp interfaces are captured remarkably well for the $\mu$ profile, even at the presence of 20% noise. The inversion for $\lambda$ is also satisfactory.

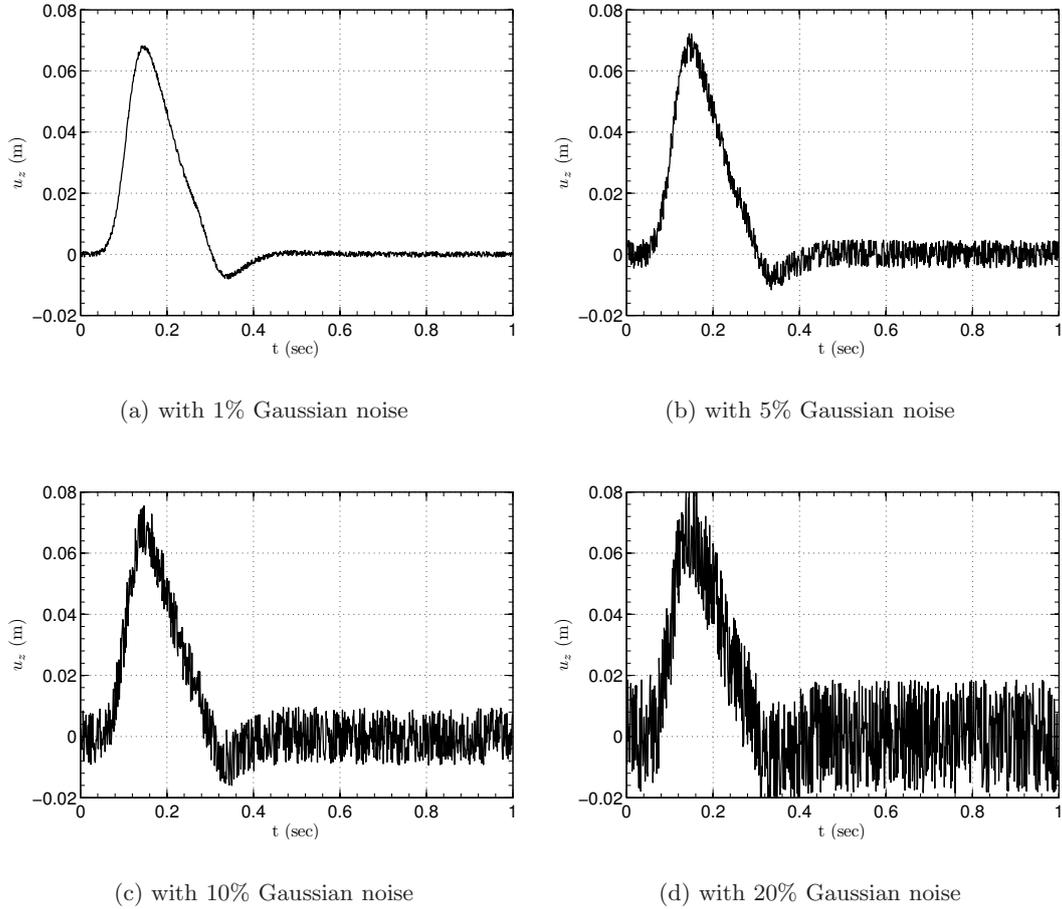

(a) with 1% Gaussian noise

(b) with 5% Gaussian noise

(c) with 10% Gaussian noise

(d) with 20% Gaussian noise

Figure 27: Measured displacement response of the layered medium with inclusion, at $(x, y, z) = (3.125, -13.75, 0)$ m, due to the Gaussian pulse $p_{20}$, contaminated with Gaussian noise.

*4.5. Layered medium with three inclusions*

In the last example, we consider a layered medium, with three inclusions, to study the performance of our inversion scheme for a more complex material profile. The problem consists of an 80 m × 80 m × 45 m medium, where a 6.25 m-thick PML is placed at its truncation boundaries. The material profiles are given by:

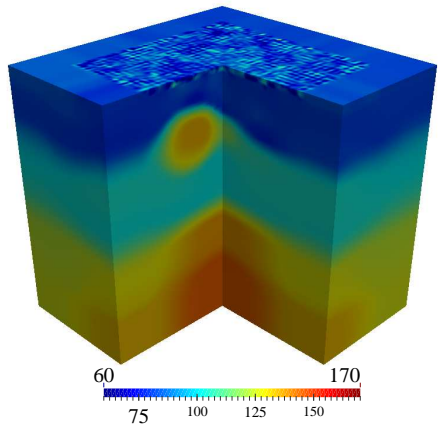
(a) $\lambda$ (1% Gaussian noise)

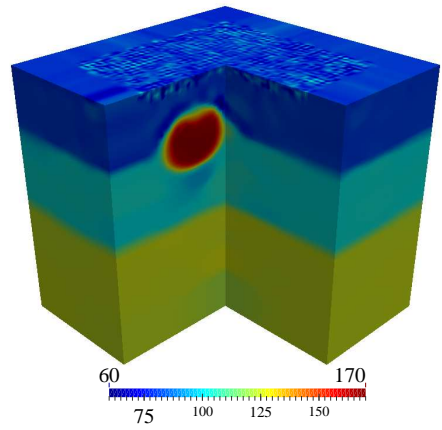
(b) $\mu$ (1% Gaussian noise)

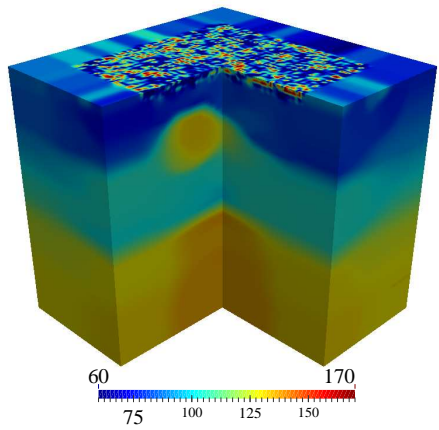
(c) $\lambda$ (5% Gaussian noise)

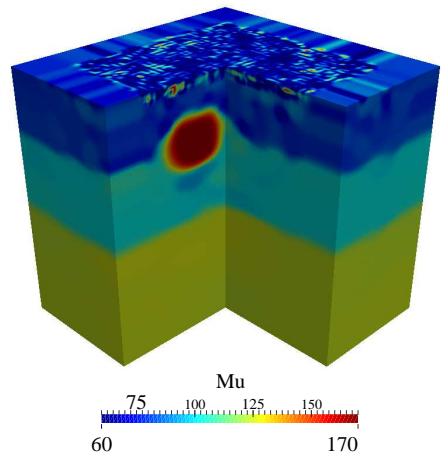
(d) $\mu$ (5% Gaussian noise)

Figure 28: Simultaneous inversion for $\lambda$ and $\mu$ using measured data contaminated with 1% and 5% Gaussian noise (layered medium with inclusion).

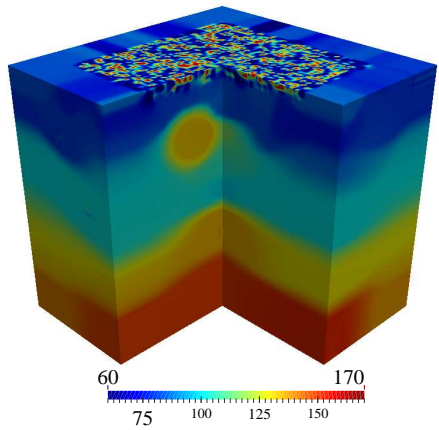

(a) $\lambda$ (10% Gaussian noise)

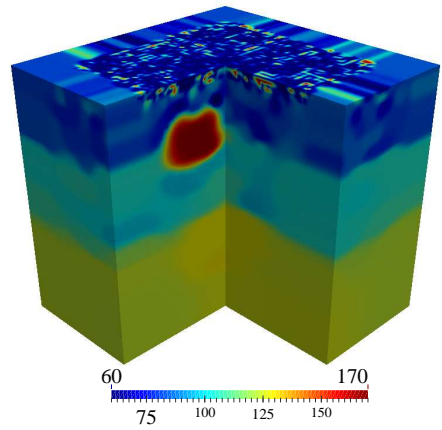

(b) $\mu$ (10% Gaussian noise)

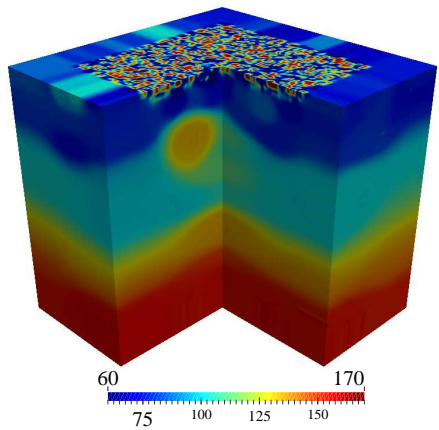

(c) $\lambda$ (20% Gaussian noise)

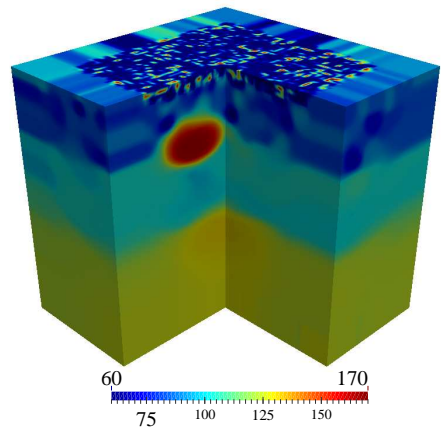

(d) $\mu$ (20% Gaussian noise)

Figure 29: Simultaneous inversion for $\lambda$ and $\mu$ using measured data contaminated with 10% and 20% Gaussian noise (layered medium with inclusion).

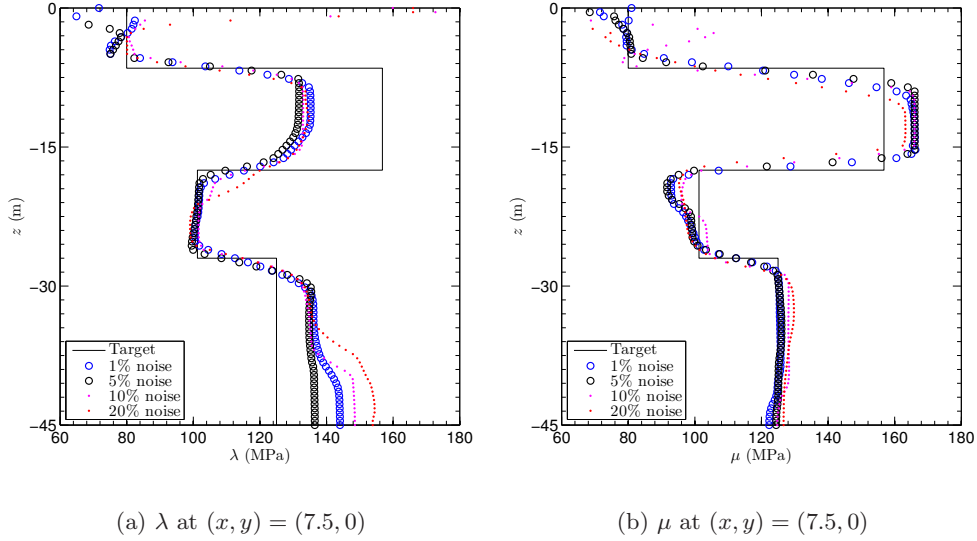

(a) $\lambda$ at $(x,y) = (7.5, 0)$        (b) $\mu$ at $(x,y) = (7.5, 0)$

Figure 30: Cross-sectional profiles of $\lambda$ and $\mu$ at different noise levels (layered medium with inclusion).

$$\lambda(z) = \mu(z) = \begin{cases} 80 \text{ MPa}, & \text{for } -15 \text{ m} \leq z \leq 0 \text{ m}, \\ 101.25 \text{ MPa}, & \text{for } -30 \text{ m} \leq z < -15 \text{ m}, \\ 125 \text{ MPa}, & \text{for } -50 \text{ m} \leq z < -30 \text{ m}, \\ 156.8 \text{ MPa}, & \text{for spheroidal: } (\frac{x+20}{3.75})^2 + (\frac{y+20}{20})^2 + (\frac{z+8.75}{3.75})^2 \leq 1, \\ 156.8 \text{ MPa}, & \text{for ellipsoidal: } (\frac{x-20}{15})^2 + (\frac{y-20}{7.5})^2 + (\frac{z+30}{5})^2 \leq 1, \\ 80 \text{ MPa}, & \text{for sphere: } (x-20)^2 + (y+20)^2 + (z+35)^2 \leq 6.25, \end{cases}$$

and are shown in Figure 31, with constant mass density $\rho = 2000$ kg/m$^3$. Figures 32(a) and 32(b) depict the target profiles on a cross-section through the domain situated at 8.75 m and 35 m from the top surface, going through the ellipsoid's and sphere's midplane, respectively. In terms of the smallest wavelength[13] the prescribed geometry comprises a domain of $16 \times 16 \times 9$ wavelengths long, wide, and deep, a spherical inclusion with a diameter 2.5 wavelengths, an ellipsoidal inclusion of $6 \times 3 \times 2$ wavelengths, and a spheroidal inclusion of $1.5 \times 8 \times 1.5$ wavelengths. The material properties at the interfaces $\Gamma^{\text{I}}$ are extended into the PML buffer. The interior and PML domains are discretized by quadratic hexahedral spectral elements (i.e., 27-noded bricks, and quadratic-quadratic pairs of approximation for displacement and stress components in the PML, and, also, quadratic approximation for material properties) of size 1.25 m, and $\Delta t = 10^{-3}$ s. This leads to $9,404,184$ state unknowns, and $2,429,586$ material parameters. To illuminate the domain, we use vertical stress loads with Gaussian pulse temporal signatures, applied on the surface of the medium over a region ($-37.5$ m $\leq x, y \leq 37.5$ m), whereas receivers are placed at every grid point, within the same region as the load.

---

[13]The smallest wavelength is equal to the smallest velocity in the formation 200 m/s, divided by the largest probing frequency 40 Hz, i.e., 5 m.

To narrow the feasibility space and alleviate difficulties with solution multiplicity, we use the Total Variation regularization, with $\epsilon = 0.01$, combined with the regularization factor continuation scheme outlined in Section 3.3.1, the source-frequency continuation scheme in Section 3.3.2, and the biasing scheme for $\lambda$ search directions in Section 3.3.3. Specifically, we use the regularization parameter $\varrho = 0.5$ when illuminating the medium with pulse $p_{20}$ for 60 iterations. Next, we use $\varrho = 0.4$ with pulse $p_{30}$ up to the $290^{\text{th}}$ iteration. Finally, we use $\varrho = 0.3$ with pulse $p_{40}$ and stop at the $741^{\text{st}}$ iteration. In all the three cases, the total simulation time is $T = 0.7$ s.

Figure 33 shows the inverted profile along a cross-section that cuts through the domain from $(x, y) = (-20, -46.5)$ to $(-20, 20)$ to $(46.5, 20)$. The cross section passes through the larger semi-principal axes of both ellipsoids, and shows very good reconstruction of the $\mu$ profile, and satisfactory inversion of the $\lambda$ profile. The layering is sharp, and the ellipsoids are captured well. In Figure 34, a cross section of the inverted profiles from $(x, y) = (20, 46.5)$ to $(20, -20)$ to $(-46.5, -20)$ is displayed, where it passes through the smaller semi-principal axes of the ellipsoids and the center of the sphere. The ellipsoids are well captured; however, the sphere, which consists of "soft" materials, can hardly be noticed, especially, in the $\lambda$ profile. Figure 35(a) and 35(b) show the inverted profiles on a cross-section through the domain, situated at 8.75 m from the surface, going through the top ellipsoid's midplane, and show satisfactory reconstruction of the ellipsoid. To see the reconstruction of the sphere in more detail, we consider a cross-section, which goes through the sphere's midplane, situated at 35 m from the top surface; this is shown in Figure 35(c) and 35(d). The sphere's footprint is visible in the $\lambda$ profile, whereas it is more conspicuous in the $\mu$ profile.

We also compare cross sections of the inverted profiles with the target along three different lines, which pass through the ellipsoids and the sphere. These are shown in Figure 36. Overall, the inverted profiles are satisfactory.

## 5. Conclusions

We discussed a full-waveform-based inversion methodology for imaging the elastic properties of a soil medium in three-dimensional, arbitrarily heterogeneous, semi-infinite domains. The problem typically arises in geotechnical site characterization and geophysical explorations, where high-fidelity imaging of the two Lamé parameters (or an equivalent pair) is of interest. Elastic waves are used as probing agents to interrogate the soil medium,

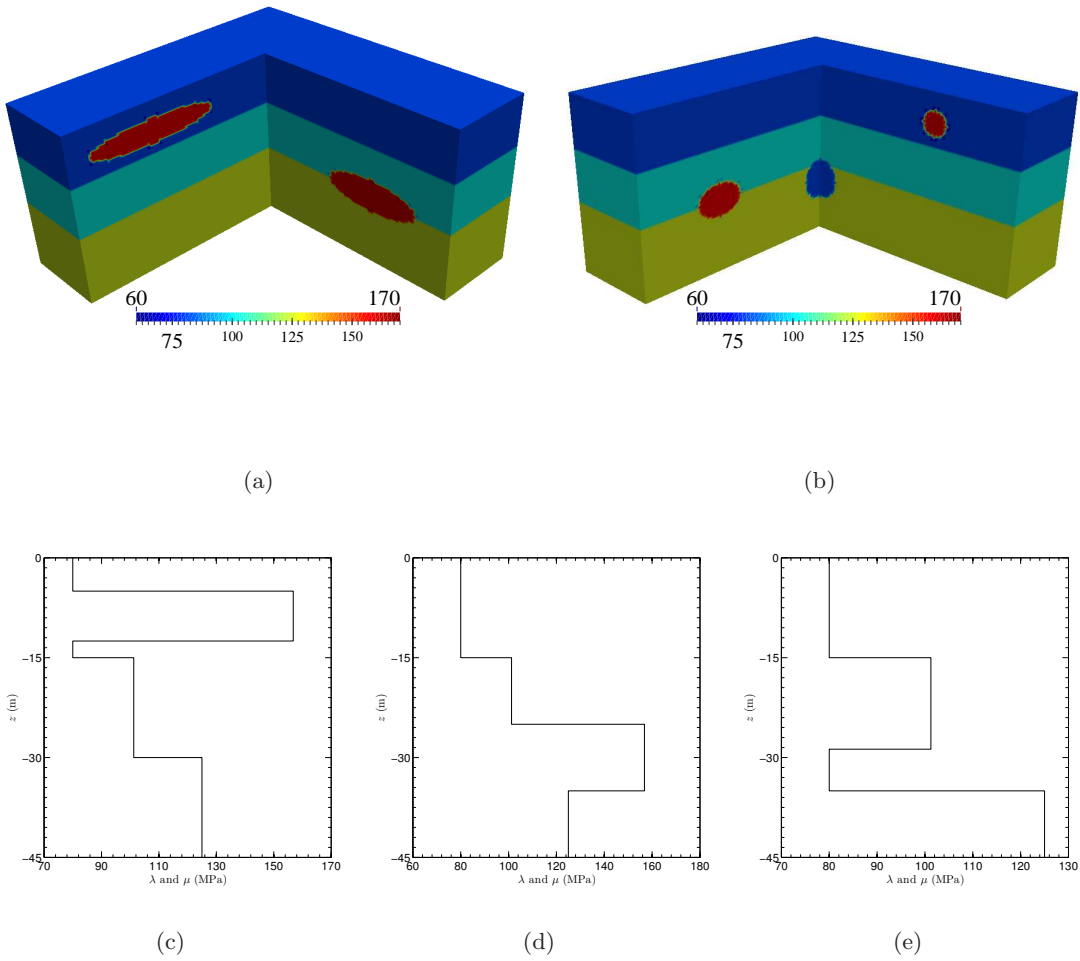

Figure 31: Layered medium with three inclusions: target $\lambda$ and $\mu$ (a) along a cross-section that cuts through the domain from $(x,y)$ = (-20,-46.5) to (-20, 20) to (46.5, 20); (b) along a cross-section that cuts through the medium from $(x,y)$ = (20, 46.5) to (20,-20) to (-46.5,-20); (c) profile at $(x,y)$ = (-20,-20); (d) profile at $(x,y)$ = (20, 20); and (e) profile at $(x,y)$ = (20,-20).

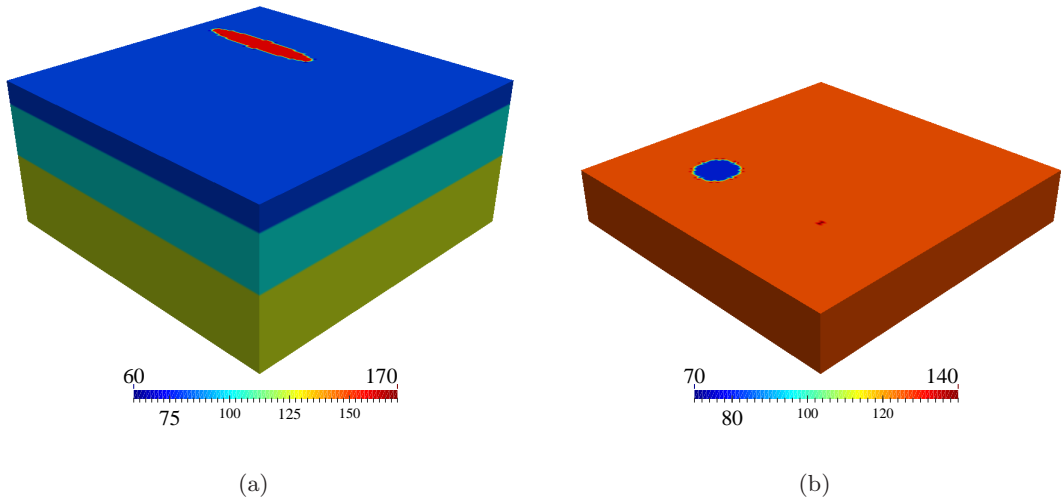

Figure 32: Layered medium with three inclusions: target $\lambda$ and $\mu$ on (a) the $z = -8.75$ m cross-section; and (b) the $z = -35$ m cross-section.

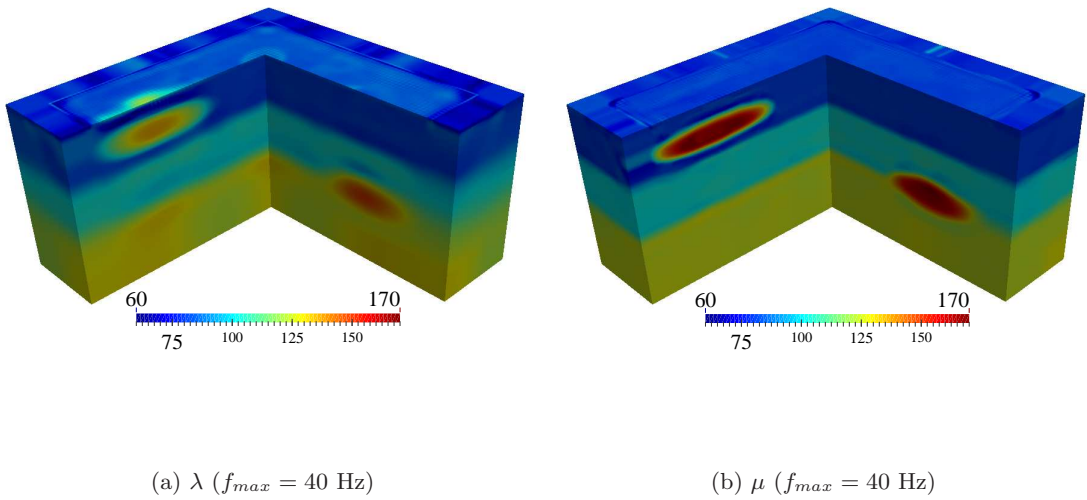

(a) $\lambda$ ($f_{max} = 40$ Hz)

(b) $\mu$ ($f_{max} = 40$ Hz)

Figure 33: Simultaneous inversion for $\lambda$ and $\mu$: cross-section cuts through the domain from $(x, y) = (-20, -46.5)$ to $(-20, 20)$ to $(46.5, 20)$ (layered medium with three inclusions).

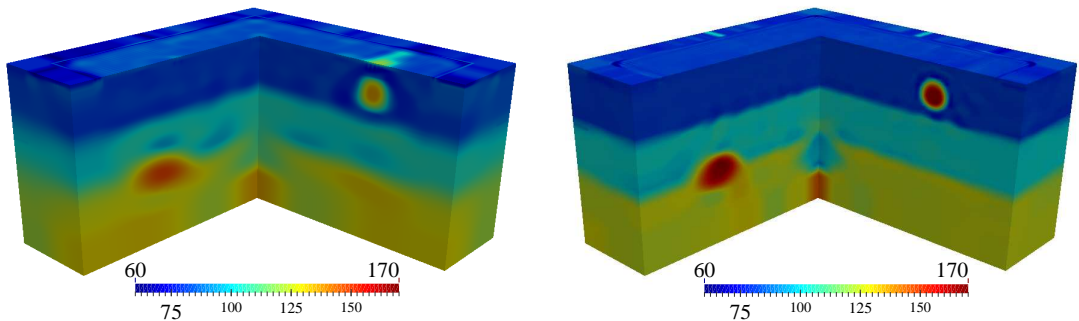

(a) $\lambda$ ($f_{max} = 40$ Hz)    (b) $\mu$ ($f_{max} = 40$ Hz)

Figure 34: Simultaneous inversion for $\lambda$ and $\mu$: cross-section cuts through the domain from $(x, y) = (20, 46.5)$ to $(20, -20)$ to $(-46.5, -20)$ (layered medium with three inclusions).

and the response of the medium to these waves are collected at receivers located on the ground surface. The inversion process relies on minimizing a misfit between the collected response at receiver locations, and a computed response based on a trial distribution of the Lamé parameters. We used the apparatus of PDE-constrained optimization to impose the forward wave propagation equations to the minimization problem, directly in the time-domain. Moreover, PMLs were used to limit the extent of the computational domain.

To alleviate the ill-posedness, associated with inverse problems, we used regularization schemes, along with a regularization factor continuation scheme, which tunes the regularization factor adaptively at each inversion iteration. We discussed additional strategies to robustify the inversion algorithm: specifically, we used (a) a source-frequency continuation scheme such that the inversion process evolves by using low-frequency sources, and, gradually, we use sources with higher frequencies; and (b) a biasing scheme for the $\lambda$-profile, such that, at early iterations of inversion, the search direction for $\lambda$ is biased based on that of $\mu$. The latter strategy, in particular, improves the reconstruction of the material profiles when simultaneous inversion of the two Lamé parameters is exercised. To the best of our knowledge, this is the first attempt that the two Lamé parameters have been successfully reconstructed in three-dimensional PML-truncated domains.

In order to resolve the forward wave propagation problem, we used a recently developed hybrid finite element approach, where a displacement-stress formulation for the PML is coupled to a standard displacement-only formulation for the interior domain, resulting in

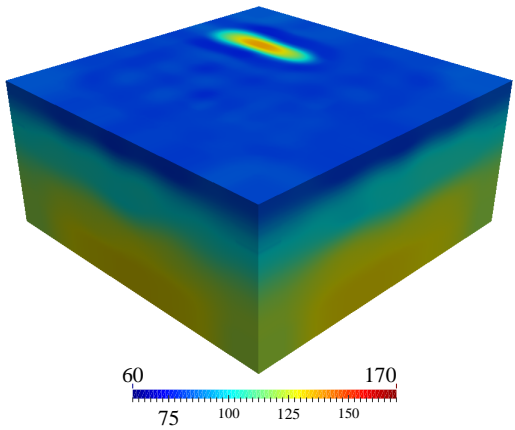
(a) $\lambda$ ($f_{max} = 40$ Hz)

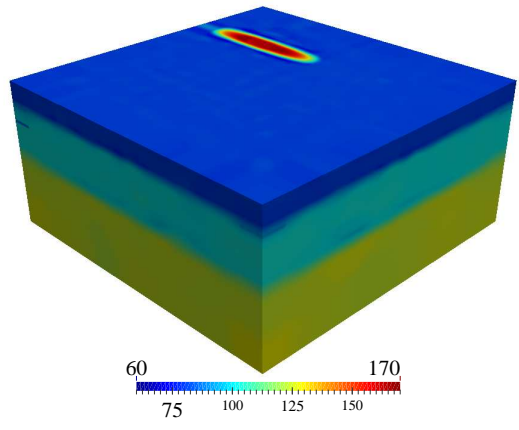
(b) $\mu$ ($f_{max} = 40$ Hz)

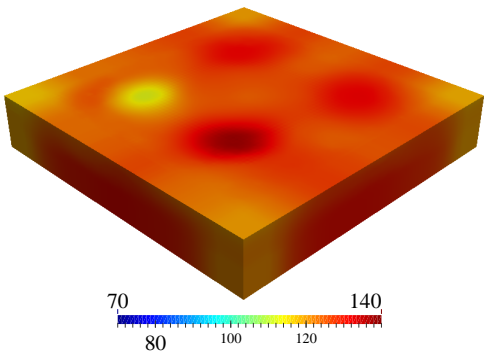
(c) $\lambda$ ($f_{max} = 40$ Hz)

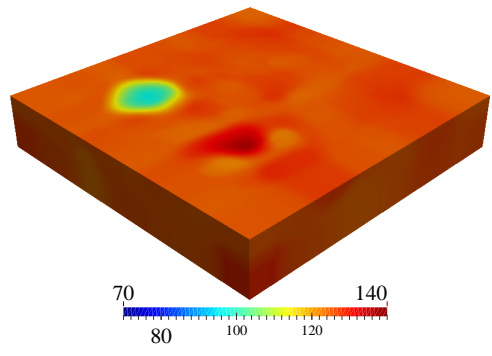
(d) $\mu$ ($f_{max} = 40$ Hz)

Figure 35: Layered medium with three inclusions: (a) inverted $\lambda$ profile on the $z = -8.75$ m cross-section; (b) inverted $\mu$ profile on the $z = -8.75$ m cross-section; (c) inverted $\lambda$ profile on the $z = -35$ m cross-section; and (d) inverted $\mu$ profile on the $z = -35$ m cross-section.

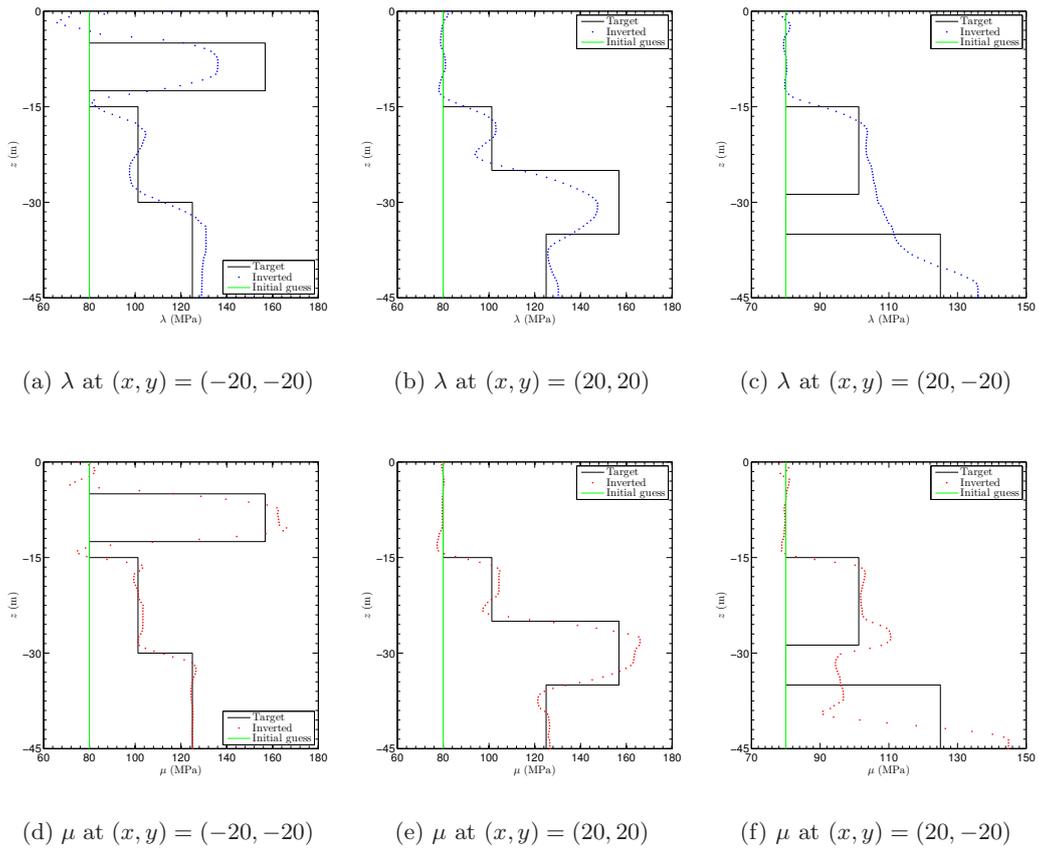

(a) $\lambda$ at $(x, y) = (-20, -20)$
(b) $\lambda$ at $(x, y) = (20, 20)$
(c) $\lambda$ at $(x, y) = (20, -20)$
(d) $\mu$ at $(x, y) = (-20, -20)$
(e) $\mu$ at $(x, y) = (20, 20)$
(f) $\mu$ at $(x, y) = (20, -20)$

Figure 36: Cross-sectional profiles of $\lambda$ and $\mu$ (layered medium with three inclusions).

a scheme with optimal computational cost. Time-integration is accomplished by using an explicit Runge-Kutta scheme, which is well-suited for large-scale problems on parallel computers.

By comparing directional finite differences of the discrete objective functional, and directional derivatives obtained via the control problems, we verified the accuracy and correctness of the material gradients. We reported numerical results demonstrating successful reconstruction of both Lamé parameters for smooth and sharp profiles. Overall, the framework discussed in this paper seems practical, and promising.

## 6. Acknowledgments

We would like to thank two anonymous reviewers whose comments improved the quality of this paper. The first author also wishes to thank Dr. Georg Stadler of Courant Institute of Mathematical Sciences for fruitful discussions and helpful comments. Partial support for the authors' research has been provided by the National Science Foundation under grant awards CMMI-0619078 and CMMI-1030728, and through an Academic Alliance Excellence grant between the King Abdullah University of Science and Technology in Saudi Arabia (KAUST) and the University of Texas at Austin. This support is gratefully acknowledged.

## Appendix A. Gradient of a functional

The gradient of a functional $\mathcal{F} : \mathcal{H} \to \mathbb{R}$, where $\mathcal{H}$ is a Hilbert space, is defined as the Riesz-representation of the derivative $\mathcal{F}'(\mathbf{q})(\tilde{\mathbf{q}})$, such that:

$$( \mathcal{G}(\mathbf{q}), \tilde{\mathbf{q}} )_\mathcal{H} = \mathcal{F}'(\mathbf{q})(\tilde{\mathbf{q}}) \qquad \forall \tilde{\mathbf{q}} \in \mathcal{H}, \tag{A.1}$$

where $\mathcal{G}$ denotes the gradient, and we use the following notation for the Gâteaux derivative of $\mathcal{F}$ at $\mathbf{q}$ in a direction $\tilde{\mathbf{q}}$:

$$\mathcal{F}'(\mathbf{q})(\tilde{\mathbf{q}}) = \lim_{h \to 0} \frac{\mathcal{F}(\mathbf{q} + h\tilde{\mathbf{q}}) - \mathcal{F}(\mathbf{q})}{h}. \tag{A.2}$$

With this definition, it is not possible to talk about the gradient, without specifying the inner-product utilized to represent the derivative [5].

## Appendix B. The adjoint problem time-integration scheme

We outline the explicit $4^{\text{th}}$-order Runge-Kutta method (RK-4) for the reverse time-integration of the adjoint problem. Upon using spectral elements, with Legendre-Gauss-Lobatto (LGL) quadrature rule, the mass-like matrix $\mathbf{M}$ becomes diagonal; therefore, its inverse can be readily computed. We use the following notation:

$$\hat{\mathbf{C}} = \mathbf{C}\,\mathbf{M}^{-1}, \qquad\qquad \hat{\mathbf{K}} = \mathbf{K}\,\mathbf{M}^{-1}, \qquad\qquad (B.1a)$$

$$\hat{\mathbf{G}} = \mathbf{G}\,\mathbf{M}^{-1}, \qquad\qquad \hat{\mathbf{f}} = \mathbf{M}^{-1}\,\mathbf{f}^{\text{adj}}. \qquad\qquad (B.1b)$$

Using (B.1), (15) becomes:

$$\frac{d}{dt}\begin{bmatrix}\mathbf{y}_0\\ \mathbf{y}_1\\ \mathbf{y}_2\end{bmatrix} = \begin{bmatrix}\mathbf{0} & \mathbf{I} & \mathbf{0}\\ \mathbf{0} & \mathbf{0} & \mathbf{I}\\ \hat{\mathbf{G}}^T & -\hat{\mathbf{K}}^T & \hat{\mathbf{C}}^T\end{bmatrix}\begin{bmatrix}\mathbf{y}_0\\ \mathbf{y}_1\\ \mathbf{y}_2\end{bmatrix} + \begin{bmatrix}\mathbf{0}\\ \mathbf{0}\\ \hat{\mathbf{f}}\end{bmatrix}. \qquad (B.2)$$

The scheme entails computing the following vectors:

$$\mathbf{k}_{10} = \mathbf{y}_1^n,$$

$$\mathbf{k}_{11} = \mathbf{y}_2^n,$$

$$\mathbf{k}_{12} = \hat{\mathbf{C}}\mathbf{y}_2^n - \hat{\mathbf{K}}\mathbf{y}_1^n + \hat{\mathbf{G}}\mathbf{y}_0^n + \hat{\mathbf{f}}^n,$$

$$\mathbf{k}_{20} = \mathbf{y}_1^n - \frac{\Delta t}{2}\,\mathbf{k}_{11},$$

$$\mathbf{k}_{21} = \mathbf{y}_2^n - \frac{\Delta t}{2}\,\mathbf{k}_{12},$$

$$\mathbf{k}_{22} = \hat{\mathbf{C}}(\mathbf{y}_2^n - \frac{\Delta t}{2}\,\mathbf{k}_{12}) - \hat{\mathbf{K}}(\mathbf{y}_1^n - \frac{\Delta t}{2}\,\mathbf{k}_{11}) + \hat{\mathbf{G}}(\mathbf{y}_0^n - \frac{\Delta t}{2}\,\mathbf{k}_{10}) + \hat{\mathbf{f}}^{n-\frac{1}{2}},$$

$$\mathbf{k}_{30} = \mathbf{y}_1^n - \frac{\Delta t}{2}\,\mathbf{k}_{21},$$

$$\mathbf{k}_{31} = \mathbf{y}_2^n - \frac{\Delta t}{2}\,\mathbf{k}_{22},$$

$$\mathbf{k}_{32} = \hat{\mathbf{C}}(\mathbf{y}_2^n - \frac{\Delta t}{2}\,\mathbf{k}_{22}) - \hat{\mathbf{K}}(\mathbf{y}_1^n - \frac{\Delta t}{2}\,\mathbf{k}_{21}) + \hat{\mathbf{G}}(\mathbf{y}_0^n - \frac{\Delta t}{2}\,\mathbf{k}_{20}) + \hat{\mathbf{f}}^{n-\frac{1}{2}},$$

$$\mathbf{k}_{40} = \mathbf{y}_1^n - \Delta t\,\mathbf{k}_{31},$$

$$\mathbf{k}_{41} = \mathbf{y}_2^n - \Delta t\,\mathbf{k}_{32},$$

$$\mathbf{k}_{42} = \hat{\mathbf{C}}(\mathbf{y}_2^n - \Delta t\,\mathbf{k}_{32}) - \hat{\mathbf{K}}(\mathbf{y}_1^n - \Delta t\,\mathbf{k}_{31}) + \hat{\mathbf{G}}(\mathbf{y}_0^n - \Delta t\,\mathbf{k}_{30}) + \hat{\mathbf{f}}^{n-1}.$$

Finally, the solution at time step $(n-1)$ can be updated via:

$$\begin{bmatrix}\mathbf{y}_0\\ \mathbf{y}_1\\ \mathbf{y}_2\end{bmatrix}^{n-1} = \begin{bmatrix}\mathbf{y}_0\\ \mathbf{y}_1\\ \mathbf{y}_2\end{bmatrix}^n - \frac{\Delta t}{6}\begin{bmatrix}\mathbf{k}_{10} + 2\,\mathbf{k}_{20} + 2\,\mathbf{k}_{30} + \mathbf{k}_{40}\\ \mathbf{k}_{11} + 2\,\mathbf{k}_{21} + 2\,\mathbf{k}_{31} + \mathbf{k}_{41}\\ \mathbf{k}_{12} + 2\,\mathbf{k}_{22} + 2\,\mathbf{k}_{32} + \mathbf{k}_{42}\end{bmatrix}. \qquad (B.4)$$

**Appendix C. Discretization of the control problems**

In Section 3.1.3, we discussed the $\lambda$- and $\mu$-control problems. In this part, we consider their spatial discretization. As discussed in [13], we use the basis $\mathbf{\Phi}$ for the spatial discretization of $\mathbf{w}(\mathbf{x}, t)$ and $\mathbf{u}(\mathbf{x}, t)$, and $\boldsymbol{\chi}$ is the basis for discretizing $\lambda(\mathbf{x})$ and $\mu(\mathbf{x})$. For instance, if we approximate $\lambda(\mathbf{x})$ with $\lambda_h(\mathbf{x})$, then $\lambda_h(\mathbf{x}) = \boldsymbol{\chi}^T \boldsymbol{\lambda}$, where $\boldsymbol{\lambda}$ comprises the vector of nodal values for $\lambda$. In the following, subscripts in the basis indicate derivatives, and $\mathbf{u}_h = (\mathbf{u}_x^T, \mathbf{u}_y^T, \mathbf{u}_z^T)^T$ and $\mathbf{w}_h = (\mathbf{w}_x^T, \mathbf{w}_y^T, \mathbf{w}_z^T)^T$ is the vector of discrete values of the state and adjoint variables, respectively. Accordingly:

$$\tilde{\mathbf{M}} = \int_{\Omega^{\text{RD}}} \boldsymbol{\chi} \boldsymbol{\chi}^T \, \mathrm{d}\Omega. \tag{C.1a}$$

For Tikhonov regularization:

$$\mathbf{g}_{\text{reg}}^{\lambda} = \int_{\Omega^{\text{RD}}} (\boldsymbol{\chi}_x \boldsymbol{\chi}_x^T + \boldsymbol{\chi}_y \boldsymbol{\chi}_y^T + \boldsymbol{\chi}_z \boldsymbol{\chi}_z^T) \boldsymbol{\lambda} \, \mathrm{d}\Omega, \tag{C.1b}$$

$$\mathbf{g}_{\text{reg}}^{\mu} = \int_{\Omega^{\text{RD}}} (\boldsymbol{\chi}_x \boldsymbol{\chi}_x^T + \boldsymbol{\chi}_y \boldsymbol{\chi}_y^T + \boldsymbol{\chi}_z \boldsymbol{\chi}_z^T) \boldsymbol{\mu} \, \mathrm{d}\Omega. \tag{C.1c}$$

For Total Variation regularization:

$$\mathbf{g}_{\text{reg}}^{\lambda} = \int_{\Omega^{\text{RD}}} \frac{(\boldsymbol{\chi}_x \boldsymbol{\chi}_x^T + \boldsymbol{\chi}_y \boldsymbol{\chi}_y^T + \boldsymbol{\chi}_z \boldsymbol{\chi}_z^T) \boldsymbol{\lambda}}{\left(\boldsymbol{\lambda}^T (\boldsymbol{\chi}_x \boldsymbol{\chi}_x^T + \boldsymbol{\chi}_y \boldsymbol{\chi}_y^T + \boldsymbol{\chi}_z \boldsymbol{\chi}_z^T) \boldsymbol{\lambda} + \epsilon\right)^{\frac{1}{2}}} \, \mathrm{d}\Omega, \tag{C.1d}$$

$$\mathbf{g}_{\text{reg}}^{\mu} = \int_{\Omega^{\text{RD}}} \frac{(\boldsymbol{\chi}_x \boldsymbol{\chi}_x^T + \boldsymbol{\chi}_y \boldsymbol{\chi}_y^T + \boldsymbol{\chi}_z \boldsymbol{\chi}_z^T) \boldsymbol{\mu}}{\left(\boldsymbol{\mu}^T (\boldsymbol{\chi}_x \boldsymbol{\chi}_x^T + \boldsymbol{\chi}_y \boldsymbol{\chi}_y^T + \boldsymbol{\chi}_z \boldsymbol{\chi}_z^T) \boldsymbol{\mu} + \epsilon\right)^{\frac{1}{2}}} \, \mathrm{d}\Omega. \tag{C.1e}$$

Moreover,

$$\mathbf{g}_{\text{mis}}^{\lambda} = -\int_0^T \int_{\Omega^{\text{RD}}} \boldsymbol{\chi} \left(\mathbf{\Phi}_x^T \mathbf{w}_x + \mathbf{\Phi}_y^T \mathbf{w}_y + \mathbf{\Phi}_z^T \mathbf{w}_z\right)\left(\mathbf{\Phi}_x^T \mathbf{u}_x + \mathbf{\Phi}_y^T \mathbf{u}_y + \mathbf{\Phi}_z^T \mathbf{u}_z\right) \mathrm{d}\Omega \, dt, \tag{C.1f}$$

$$\mathbf{g}_{\text{mis}}^{\mu} = -\int_0^T \int_{\Omega^{\text{RD}}} \boldsymbol{\chi} \left( 2 \left(\mathbf{\Phi}_x^T \mathbf{w}_x \, \mathbf{\Phi}_x^T \mathbf{u}_x + \mathbf{\Phi}_y^T \mathbf{w}_y \, \mathbf{\Phi}_y^T \mathbf{u}_y + \mathbf{\Phi}_z^T \mathbf{w}_z \, \mathbf{\Phi}_z^T \mathbf{u}_z\right) \right.$$
$$+ (\mathbf{\Phi}_y^T \mathbf{w}_x + \mathbf{\Phi}_x^T \mathbf{w}_y)(\mathbf{\Phi}_x^T \mathbf{u}_y + \mathbf{\Phi}_y^T \mathbf{u}_x) + (\mathbf{\Phi}_z^T \mathbf{w}_x + \mathbf{\Phi}_x^T \mathbf{w}_z)(\mathbf{\Phi}_x^T \mathbf{u}_z + \mathbf{\Phi}_z^T \mathbf{u}_x)$$
$$\left. + (\mathbf{\Phi}_z^T \mathbf{w}_y + \mathbf{\Phi}_y^T \mathbf{w}_z)(\mathbf{\Phi}_y^T \mathbf{u}_z + \mathbf{\Phi}_z^T \mathbf{u}_y) \right) \mathrm{d}\Omega \, dt. \tag{C.1g}$$

In (C.1a), upon using spectral elements with LGL quadrature rule, $\tilde{\mathbf{M}}$ becomes diagonal; thus, its inverse can be computed easily.